\RequirePackage{ifthen}
\newboolean{MPA}
\setboolean{MPA}{false}

\ifthenelse {\boolean{MPA}}
{
\documentclass{svjour3}
\smartqed
\usepackage[margin=1.3in]{geometry}
} {
\documentclass[11pt]{article}
\usepackage[margin=1.0in]{geometry}
}

\usepackage{graphicx}
\usepackage{amssymb}
\usepackage{amsmath}
\usepackage{verbatim,booktabs}
\usepackage{subfigure,epsfig,url,psfrag}
\usepackage{float}
\usepackage{mathrsfs}

\usepackage[usenames,dvipsnames]{pstricks}
\usepackage{epsfig}

\usepackage{cite} 

\usepackage{color}
\usepackage{enumerate}
\usepackage{enumitem}

\usepackage{todonotes}

\usepackage{hyperref}
\usepackage[capitalize,noabbrev]{cleveref}

\crefname{question}{Question}{Questions}
\crefname{step}{Step}{Steps}
\crefname{claim}{Claim}{Claims}
\crefname{problem}{Problem}{Problems}
\crefname{definition}{Definition}{Definitions}
\crefname{observation}{Observation}{Observations}

\newcommand{\pare}[1]{\left(#1\right)}
\newcommand{\bra}[1]{\left\{#1\right\}}
\newcommand{\card}[1]{\left|#1\right|}

\newcommand{\ie}{i.e., }

\newcommand{\NP}{{\mathcal NP}}
\newcommand{\ME}{E_{\max}}
\newcommand{\R}{{\mathbb R}}

\renewcommand{\P}{\mathcal P}
\renewcommand{\S}{\mathcal S}

\def\01{\ensuremath{0\mathord{-}1}}
\allowdisplaybreaks

\DeclareMathOperator{\conv}{conv}
\DeclareMathOperator{\poly}{poly}
\DeclareMathOperator{\proj}{proj}

\DeclareMathOperator{\cl}{cl}
\DeclareMathOperator{\MP}{MP}
\DeclareMathOperator{\MPLP}{MP^{LP}}
\DeclareMathOperator{\BQP}{BQP}
\DeclareMathOperator{\BQPLP}{BQP^{LP}}
\DeclareMathOperator{\CER}{CER}

\ifthenelse {\boolean{MPA}}
{

\newenvironment{prf}[1][]
{\begin{proof}}
{\qed \end{proof}}

\newenvironment{prfc}[1][]
{\begin{proof}[#1]}
{\qed \end{proof}}

\newenvironment{prfh}[1][]
{\begin{proof}}
{\end{proof}}

\newcounter{claim} 
\renewenvironment{claim}[1][]
{\refstepcounter{claim} \begin{trivlist} \item[] {\bf Claim~\theclaim}\space#1 \itshape}
{\end{trivlist}}

\newenvironment{cpf}
{\begin{trivlist} \item[] {\em Proof of claim }}
{$\hfill\diamond$ \end{trivlist}}

\journalname{Mathematical Programming A}

\newtheorem{observation}{Observation}

} {

\usepackage{amsthm}
\newtheorem{theorem}{Theorem}
\newtheorem{definition}{Definition}
\newtheorem{lemma}{Lemma}

\newtheorem{proposition}{Proposition}
\newtheorem{example}{Example}
\newtheorem{observation}{Observation}
\newtheorem{claim}{Claim}

\newenvironment{prf}[1][]
{\begin{proof}}
{\end{proof}}

\newenvironment{cpf}
{\begin{trivlist} \item[] {\em Proof of claim. }}
{$\hfill\diamond$ \end{trivlist}}

}

 \allowdisplaybreaks

\date{June 10, 2026}

\title{The complete edge relaxation for binary polynomial optimization}

\ifthenelse {\boolean{MPA}}
{
\titlerunning{The complete edge relaxation for binary polynomial optimization}

\author{Alberto Del Pia \and Aida Khajavirad}

\institute{Alberto Del Pia \at
              Department of Industrial and Systems Engineering \& Wisconsin Institute for Discovery, 
              University of Wisconsin-Madison.
              E-mail: {\tt delpia@wisc.edu}.
           \and
           Aida Khajavirad \at
              Department of Industrial and Systems Engineering,
              Lehigh University.
              E-mail: {\tt aida@lehigh.edu}.
}
}
{
\author{Alberto Del Pia
\thanks{Department of Industrial and Systems Engineering \& Wisconsin Institute for Discovery,
             University of Wisconsin-Madison.
             E-mail: {\tt delpia@wisc.edu}.
             }
\and
Aida Khajavirad
\thanks{Department of Industrial and Systems Engineering,
             Lehigh University.
             E-mail: {\tt aida@lehigh.edu}.
             }
}
}

\begin{document}

\maketitle


\begin{abstract}
    We consider the multilinear polytope, defined as the convex hull of the feasible region of a lifted binary polynomial optimization problem. We define a relaxation in an extended space for this polytope, which we call the complete edge relaxation. The complete edge relaxation is stronger than several well-known relaxations of the multilinear polytope, including the standard linearization, the flower relaxation, and the intersection of all possible recursive McCormick relaxations. In addition, for fixed-degree binary polynomial optimization problems, the case of primary practical interest, the complete edge relaxation is of polynomial size and is computationally efficient in practice.
    We prove that the complete edge relaxation is an extension of the multilinear polytope if and only if the corresponding hypergraph is $\alpha$-acyclic, the most general type of hypergraph acyclicity. This is in stark contrast with the widely-used standard linearization, which describes the multilinear polytope if and only if the hypergraph is Berge-acyclic, the most restrictive type of hypergraph acyclicity. Finally, we introduce a new class of facet-defining inequalities for the multilinear polytope of $\alpha$-cycles of length three, which serve as the generalization of the well-known triangle inequalities for the Boolean quadric polytope.
\ifthenelse {\boolean{MPA}}
{
\keywords{Binary polynomial optimization \and Multilinear polytope \and Complete edge relaxation \and Hypergraph acyclicity \and Extended formulations \and Generalized triangle inequalities}
\subclass{MSC 90C09 \and 90C10 \and 90C26 \and 90C57}
} {}
\end{abstract}

\medskip
\ifthenelse {\boolean{MPA}}
{}{
\emph{Key words:} Binary polynomial optimization; Multilinear polytope; Complete edge relaxation; Hypergraph acyclicity; Extended formulations; Generalized triangle inequalities.
}

\section{Introduction}
\label{sec intro}


\paragraph{Problem statement.} Binary polynomial optimization, \ie the problem of finding a binary vector maximizing a multivariate polynomial function, is a fundamental $\NP$-hard problem in discrete optimization~\cite{BorHam02,Pun22book}. To formally define this problem, we employ a hypergraph representation scheme that we first introduced in~\cite{dPKha17MOR}. A \emph{hypergraph} $G$ is a pair $(V,E)$, where $V$ is a finite set of nodes and $E$ is a set of subsets of $V$ of cardinality at least two, called the edges of $G$. 
With any hypergraph $G= (V,E)$, and the cost vector $c \in \R^{V \cup E}$, we associate the~\emph{binary polynomial optimization} problem:
\begin{align}
\label[problem]{prob BPO}
\tag{BPO}
\begin{split}
\max_{z \in \{0,1\}^V} & \quad \sum_{v\in V} {c_v z_v} + \sum_{e\in E} {c_e \prod_{v\in e} {z_v}}, 
\end{split}
\end{align}
where, without loss of generality, we assume that each node is contained in at least one edge, and $c_e \neq 0$ for all $e \in E$.
With the objective of constructing LP relaxations for~\cref{prob BPO}, we proceed with linearizing the objective function by introducing a new variable for each product term to obtain an equivalent reformulation of \cref{prob BPO} in a lifted space:
\begin{align}
\label[problem]{prob LBP}
\tag{$\ell$-BPO}
\begin{split}
\max & \qquad \sum_{v\in V} {c_v z_v} + \sum_{e\in E} {c_e z_e} \\
\text{s.t.} & \qquad z_e = \prod_{v\in e} {z_v} \qquad \forall e \in E, \\  
& \qquad z_v \in \{0,1\} \qquad \forall v \in V.
\end{split}
\end{align}
To understand the complexity of~\cref{prob BPO}, we study the facial structure of the convex hull of its feasible region. In the same vein as our previous work~\cite{dPKha17MOR}, we define the~\emph{multilinear set} as the feasible region of \cref{prob LBP} and refer to its convex hull as the~\emph{multilinear polytope}:
\begin{equation*} 
 \MP(G):= \conv\Big\{ z \in \{0,1\}^{V \cup E} : z_e = \prod_{v \in e} {z_{v}}, \; \forall e \in E \Big\}.
 \end{equation*}

\paragraph{The McCormick relaxation.} If the objective function of \cref{prob BPO} is quadratic, then the hypergraph $G$ is a graph, and the multilinear polytope coincides with the well-known \emph{Boolean quadric polytope} $\BQP(G)$ introduced first by Padberg in 1989 in the context of binary quadratic optimization~\cite{Pad89}, and later investigated by many others to study the maximum cut problem~\cite{BarMah86,DezLau97}. In~\cite{Pad89}, Padberg considered a simple relaxation of $\BQP(G)$ obtained by replacing each bilinear equation
$z_{\{u,v\}} = z_u z_v$ over $\{0,1\}^2$ by its convex hull: 
\begin{align*}
    \BQPLP(G) =
    \Big\{z \in \R^{V \cup E}:\; & z_{\{u,v\}} \geq 0, \; z_{\{u,v\}} \geq z_u + z_v -1, \; z_{\{u,v\}} \leq z_u, \; z_{\{u,v\}} \leq z_v, \\
    &\forall \{u,v\} \in E\Big\}.
\end{align*}
The polytope $\BQPLP(G)$ is often referred to as the \emph{McCormick relaxation} of $\BQP(G)$. We should remark that the inequalities defining  $\BQPLP(G)$ have been considered in earlier studies; see for example~\cite{GloWol74}.
Padberg then proved the following result regarding the tightness of the McCormick relaxation:
\begin{proposition}[proposition~8 in~\cite{Pad89}]
    $\BQP(G) = \BQPLP(G)$ if and only if $G$ is an acyclic graph.
\end{proposition}
Observe that $\BQPLP(G)$ is defined by $4|E|$ inequalities, implying that binary quadratic optimization over acyclic graphs can be solved in polynomial-time. 
Padberg then considered the simplest graph for which $\BQP(G) \subsetneq \BQPLP(G)$, \ie a cycle of length three, and introduced \emph{triangle inequalities}, a class of facet-defining inequalities for the Boolean quadric polytope whose addition significantly strengthens $\BQPLP(G)$.
In~\cite{BorCraHam92} the authors show that triangle inequalities are Chv\'atal-Gomory (CG) cuts for the McCormick relaxation. Furthermore, if $G$ is a complete graph, they proved that the addition of triangle inequalities to the McCormick relaxation gives the CG-closure of $\BQPLP(G)$. 


\paragraph{The standard linearization.} It is natural to ask whether the multilinear polytope of acyclic hypergraphs has a simple structure as well.
Interestingly, unlike graphs, the notion of cycle and acyclicity in hypergraphs is not unique. 
The most well-known types of acyclic hypergraphs, in increasing order of generality, are Berge-acyclic, $\gamma$-acyclic, $\beta$-acyclic, and $\alpha$-acyclic hypergraphs~\cite{fagin83,BeeFagMaiYan83,Dur12}. 
In a similar vein to the McCormick relaxation, a simple relaxation of the multilinear polytope can be constructed by replacing each multilinear equation $z_e = \prod_{v\in e}{z_v}$ over $z_v \in \{0,1\}$ for all $v \in e$ by its convex hull:
\begin{align*}
\MPLP(G) =\Big\{z \in \R^{V \cup E}: \;
& z_v \leq 1, \; \forall v \in V, \; 
z_e \geq 0, \; z_e \geq \sum_{v\in e}{z_v}-|e|+1, \;\forall e \in E, \\
& z_e \leq z_v, \; \forall e \in E, \; \forall v \in e\Big\}.
\end{align*}
The above relaxation is known as the \emph{standard linearization} and has been widely used in the literature~\cite{GloWol74,yc93}. 
The standard linearization has often been regarded as the generalization of McCormick relaxation for higher-degree binary polynomial optimization problems.  In~\cite{dPKha18SIOPT}, we prove the following result regarding the tightness of the standard linearization:
\begin{theorem}[theorem~7 in~\cite{dPKha18SIOPT}]
\label{th gamma}
    $\MP(G) = \MPLP(G)$ if and only if $G$ is a Berge-acyclic hypergraph.
\end{theorem}
Henceforth, given a hypergraph $G$, we define the \emph{rank} of $G$, denoted by $r$, as the maximum cardinality of an edge in $E$. Observe that $\MPLP(G)$ is defined by at most $|V|+ (r+2) |E|$ inequalities, implying that~\cref{prob BPO} over Berge-acyclic hypergraphs can be solved in polynomial time. 
However, recall that Berge-acyclicity is the most restrictive type of hypergraph acyclicity, and indeed, the standard linearization is a very weak relaxation for the multilinear polytope of general hypergraphs.
Next, we introduce a different generalization of the McCormick relaxation that is significantly stronger than the standard linearization.



\paragraph{The complete edge relaxation.}
A hypergraph $G$ is \emph{complete} if its edge set consists of all possible subsets of $V$ of cardinality at least two.
Given a hypergraph $G=(V,E)$, we say that an edge $e \in E$ is \emph{maximal} for $G$ if it is not contained in any other edge in $E$. 
Clearly, in a graph, every edge is a maximal edge. 
Then, for a graph $G$, the McCormick relaxation of the Boolean quadric polytope $\BQP(G)$ can be described as follows. 
For each maximal edge $e$ of $G$, consider the complete graph with the node set $e$, denoted by $K_e$. The McComick relaxation is then obtained by putting together the descriptions of $\BQP(K_e)$ for all $e \in E$.
While this interpretation of the McCormick relaxation seems unnecessarily complicated, applying the same scheme to hypergraphs of rank $r > 2$ yields an interesting relaxation of the multilinear polytope, in a higher-dimensional space, which we describe next.

To work with relaxations in higher-dimensional spaces, we need to introduce some terminology.
An \emph{extension} of a polyhedron $P \subseteq \R^n$ is a polyhedron $Q \subseteq \R^{n+d}$, with $d \ge 0$, such that 
$P = \bra{x \in Q : \exists y \in \R^d, \ (x,y) \in Q}$.
Any description of $Q$ by linear equalities and inequalities is called an \emph{extended formulation} of $P$. 
In this paper, the \emph{size} of an extended formulation is the size of its defining system of linear inequalities, as detailed in \cite{SchBookIP}. 

Let $G = (V,E)$ be a hypergraph and denote by $\ME \subseteq E$ the set of maximal edges of $G$. For each $e \in \ME$, denote by $K_e$ the complete hypergraph with the node set $e$. 
We can construct a relaxation for $\MP(G)$, in a higher-dimensional space, by putting together the descriptions of $\MP(K_e)$ for all $e \in \ME$. The explicit description of $\MP(K_e)$ can be readily obtained using the Reformulation Linearization Technique (RLT)~\cite{SheAda90}.
To formally define this relaxation, we next present the explicit description for the multilinear polytope of a complete hypergraph. 

\begin{proposition}[theorem~2 in~\cite{SheAda90}]\label{prop:RLT}  
Let $G=(V,E)$ be a complete hypergraph. For any $U \subseteq V$, define
\begin{equation}\label{defRLT}
    \psi(U,V):=\sum_{\substack{W \subseteq U: \\ |W| \text{ even}}}{z_{(V\setminus U)\cup W}}-\sum_{\substack{W \subseteq U: \\ |W| \text{ odd}}}{z_{(V\setminus U)\cup W}},
\end{equation}
where we define $z_{\emptyset} := 1$. Then the multilinear polytope $\MP(G)$ is given by
\begin{equation}\label{eq:rlt}
\psi(U, V) \geq 0, \quad \forall U \subseteq V.
\end{equation}
\end{proposition}
Henceforth, given a hypergraph $G=(V,E)$, we define the \emph{completion of $G$} as the hypergraph $\cl(G)=(V, \cl(E))$, where $\cl(E):= \{f\subseteq e: \; |f| \geq 2, \; \forall e \in E\}$.  
We then define the \emph{complete edge relaxation} of $\MP(G)$, denoted by $\CER(G)$, as the polytope defined by the following system of inequalities:
\begin{align*}
    \CER(G) = \Big\{z \in \R^{V \cup \cl(E)}: \psi(U, e) \geq 0, \; \forall U \subseteq e, \; \forall e \in E_{\max}\Big\},
\end{align*}
where $\psi(U,e)$ is defined by~\eqref{defRLT}. 
Notice that we have $\CER(G) = \CER(\cl(G))$.
In this paper, we settle the tightness of the complete edge relaxation. Our main result is as follows:

\begin{theorem}
\label{th main}
$\CER(G)$ is an extension of $\MP(G)$ if and only if $G$ is an $\alpha$-acyclic hypergraph.
\end{theorem}

Therefore, the complete edge relaxation is tight if and only if the hypergraph $G$ exhibits the most general form of acyclicity.
From~\cref{prop:RLT}, it follows that $\CER(G)$ is defined by at most $2^r|\ME|$ variables and inequalities. 
Therefore, if the rank of the hypergraph is bounded; \ie $r \in O(\log\poly(|V|,|E|))$, where $\poly(|V|,|E|)$ denotes a polynomial function in $|V|,|E|$, then the complete edge relaxation is of polynomial-size and can be constructed in polynomial time. This implies that~\cref{prob BPO} over $\alpha$-acyclic hypergraphs of bounded rank can be solved in polynomial-time.  
Recall that the rank of the hypergraph $G$ corresponds to the degree of the binary polynomial optimization problem, which is a small number in almost all applications; that is, in almost all application, we have $|V| \gg r$ and $|E| \gg r$. Hence, assuming $r \in O(\log\poly(|V|,|E|))$ is reasonable for all practical purposes. 
In~\cite{dPDiG23ALG}, the authors prove that solving~\cref{prob BPO} over $\alpha$-acyclic hypergraphs is, in general, strongly $\NP$-hard. This result implies that, unless $\P = \NP$, one cannot construct, in polynomial time, a polynomial-size extended formulation of the multilinear polytope of $\alpha$-acyclic hypergraphs. 

The `if' statement in~\cref{th main} is the simpler one.
In fact, in \cite{Lau09,BieMun18}, the authors present a similar result but stated in terms of the concept of \emph{primal treewidth}, rather than $\alpha$-acyclicity.
Moreover, the connection between treewidth and $\alpha$-acyclicity is detailed in \cite{dPKha21MOR}, which together with the result of~\cite{Lau09,BieMun18} implies the `if' statement in~\cref{th main}.
In this paper, we present a simple proof of this result that does not involve treewidth and makes direct use of $\alpha$-acyclicity.
Our main contribution, however, is to prove the `only if' direction; 
that is, to show that $\alpha$-acyclicity is \emph{necessary} for the tightness of the complete edge relaxation. Proving this result is fairly involved and requires the building of new tools in both hypergraph theory and polyhedral theory. First, we obtain a novel characterization of $\alpha$-acyclic hypergraphs in terms of \emph{simple cycles}. This new characterization is essential to prove our main result and is of independent interest. Second, we introduce several operations on hypergraphs, namely node fixing, node contraction, and node expansion; we
show how the multilinear polytope and the complete edge relaxation are transformed under these operations. Through a recursive application of these operations, we are able to show that it is enough to limit our attention to hypergraphs with simple cycles of length three, a structure that can be tackled directly.

In this paper, we propose the complete edge relaxation
as \emph{the generalization} of McCormick relaxation for binary polynomial optimization problems of bounded degree. While~\cref{th main} serves as the main theoretical justification for the strength of the complete edge relaxation, this relaxation has several other desirable properties, which we briefly describe next. 
In~\cite{KhaWang25} the authors prove that the existing cutting planes for the multilinear polytope, \ie flower inequalities~\cite{dPKha18SIOPT}, and running-intersection inequalities~\cite{dPKha21MOR} are implied by the complete edge relaxation. By the results in~\cite{dPKha21MOR,Kha22,SchWal24}, this further implies that, for general hypergraphs, the complete edge relaxation is stronger than the intersection of all possible \emph{recursive McCormick relaxations} of $\MP(G)$. 
On the computational side, the authors of~\cite{KhaWang25} show that for degree three, four, and five binary polynomial optimization problems arising from image restoration and decoding error-correcting codes, optimizing over $\CER(G)$ leads to tightest bounds and can be done very efficiently. 
Interestingly, for the image restoration problem, their computational experiments indicate that even though the corresponding hypergraphs are not $\alpha$-acyclic, optimizing over the complete edge relaxation almost always yields a binary solution.

Next, we show that the complete edge relaxation can further be used to obtain new classes of facet-defining inequalities for the multilinear polytope of hypergraphs that do not exhibit any degree of acyclicity.
Namely, we introduce a new class of facet-defining inequalities for the multilinear polytope of hypergraphs with $\alpha$-cycles of length three. These inequalities serve as the generalization of the well-known triangle inequalities~\cite{Pad89}. We prove that the proposed inequalities are CG-cuts for the complete edge relaxation, while they are not CG-cuts for the standard linearization, hence, demonstrating the superiority of the complete edge relaxation.


\paragraph{Organization.} The remainder of this paper is organized as follows.
In \cref{sec characterization}, we introduce a novel characterization of $\alpha$-acyclic hypergraphs, which is key to prove the reverse direction in~\cref{th main}.
Sections~\ref{sec sufficient}-\ref{sec necessary} are devoted to proving~\cref{th main}.
Namely, in~\cref{sec sufficient}, we prove the sufficiency of $\alpha$-acyclicity, while 
in~\cref{sec necessary}, we prove the necessity of $\alpha$-acyclicity. In~\cref{sec: inequalities} we introduce a new class of facet-defining inequalities for the multilinear polytope corresponding to $\alpha$-cycles of length three.

\section{A new characterization of $\alpha$-acyclic hypergraphs}
\label{sec characterization}

As we mentioned earlier, $\alpha$-acyclic hypergraphs are the most general type of acyclic hypergraphs.
In~\cite{fagin83}, the author presents twelve equivalent characterizations of $\alpha$-acyclic hypergraphs, which have been used widely ever since. We start by presenting one of these characterizations that we will use to prove the sufficiency of $\alpha$-acyclicity in~\cref{th main}:

\begin{definition}[$\alpha$-acyclic hypergraph]
    Let $G$ be a hypergraph. We say that $G$ is $\alpha$-acyclic if its edge set has the running-intersection property.
\end{definition}
Recall that a set $F$ of subsets of a finite set $V$ has the \emph{running intersection property} if there exists an ordering $f_1, f_2, \ldots, f_m$ of the sets in $F$ such that
for each $k=2, \ldots, m$, there exists $j<k$ such that 
$$f_k \cap\Big(\bigcup_{i<k} f_i\Big) \subseteq f_j.$$
Such an ordering is called a \emph{running intersection ordering} of $F$.  
Interestingly, none of the characterizations in~\cite{fagin83} is in terms of $\alpha$-cycles. More than twenty years later, in~\cite{wang05}, for the first time, the author introduced the notion of $\alpha$-cycles and proved that a hypergraph $G$ is $\alpha$-acyclic if and only if it does not contain any $\alpha$-cycles.
Below, we give the definition of $\alpha$-cycles as defined in~\cite{wang05}.


\begin{definition}[$\alpha$-cycle]
\label{def alpha}
Let $G=(V,E)$ be a hypergraph.
Let
$C= e_1, e_2, \ldots, e_{\ell}, e_{\ell+1}$
 be a sequence of edges with $e_1 = e_{\ell+1}$ for some $\ell \ge 3$.
Define $s_i := e_i \cap e_{i+1}$ for all $i \in [\ell]$.
We say that $C$ is an \emph{$\alpha$-cycle} of $G$ of length $\ell$ if the following conditions are satisfied:
\begin{align}
\tag{A1}
\label{eq alpha diff}
& e_i \setminus e_j \neq \emptyset && \quad \forall i\neq j \in [\ell],\\
\tag{A2}
\label{eq alpha cond}
& (s_{i-1}\cup s_i \cup s_{i+1}) \setminus e \neq \emptyset && \quad \forall i \in [\ell], \; \forall e \in E,
\end{align}
where the indices in \eqref{eq alpha cond} are modulus $\ell$.
\end{definition}

\begin{observation}
\label{obs simple}
Consider the sets $s_i$, $i \in [\ell]$ in~\cref{def alpha}.
We have $s_i \neq \emptyset$ for all $i \in [\ell]$.
\end{observation}

\begin{prf}
For ease of notation, in this proof the indices are modulus $\ell$.
Condition \eqref{eq alpha cond} implies
\begin{equation*}
(s_i \cup s_{i+1} \cup s_{i+2}) \setminus e_{i+2} \neq \emptyset.
\end{equation*}
However, $s_{i+1} \cup s_{i+2} \subseteq e_{i+2},$
thus $s_i \setminus e_{i+2} \neq \emptyset$, implying that $s_i \neq \emptyset$.
\end{prf}



\begin{figure}
\centering
\includegraphics[scale=0.6]{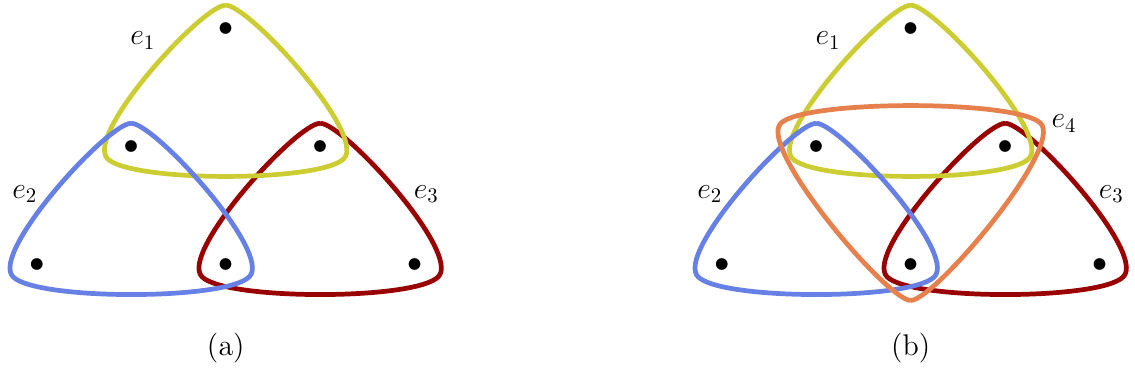}
\caption{Removing $\alpha$-cycles by adding edges to the hypergraph. (a) The hypergraph has an $\alpha$-cycle $e_1,e_2,e_3,e_1$. (b) The hypergraph is $\alpha$-acyclic because $e_4 \supseteq s_1 \cup s_2 \cup s_3$.}
\label{fig1}
\end{figure}

From~\cref{def alpha} one observes a stark contrast between $\alpha$-cycles and other types of cycles in hypergraphs, \ie Berge-cycles, $\gamma$-cycles, and $\beta$-cycles. 
Namely, while all other types of cycles are defined in terms of the edges of the cycle, the definition of an $\alpha$-cycle depends on the \emph{entire} edge set of the hypergraph. 
That is, while all other cycles are defined ``locally'', $\alpha$-cycles are defined ``globally''; \ie over the entire hypergraph. This implies that $\alpha$-cycles have a significantly more complex structure than other types of cycles. Additionally, notice that given a graph (resp. hypergraph) $G$, one can obtain an acyclic graph (resp, Berge-acyclic, $\gamma$-acyclic, $\beta$-acyclic hypergraph) by removing edges from $G$ and hence removing cycles. In contrast, to remove $\alpha$-cycles and obtain an $\alpha$-acyclic hypergraph, one could \emph{add} edges to $G$. For example, from~\eqref{eq alpha cond} it follows that if $G$ has an edge containing all its node set, then $G$ has no $\alpha$-cycles. See Figure~\ref{fig1} for an illustration.

In this paper, we present a novel characterization of $\alpha$-acyclicity, based on a ``simpler'' notion of cycles. This new notion will enable us to prove the necessity of $\alpha$-acyclicity in~\cref{th main}. We start by defining a new type of cycle in hypergraphs, which we refer to as \emph{simple cycle}.

\begin{definition}[simple-cycle]
\label{def simple}
Let $G=(V,E)$ be a hypergraph.
Let $C= e_1, e_2, \ldots, e_{\ell}, e_{\ell+1}$
be a sequence of edges with $e_1 = e_{\ell+1}$ for some $\ell \ge 3$.
Define $s_i := e_i \cap e_{i+1}$ for all $i \in [\ell]$.
We say that $C$ is a \emph{simple-cycle} of $G$, and of length $\ell$ if:
\begin{equation}
\tag{S}
\label{eq simple}
(s_i \cup s_j \cup s_k) \setminus e \neq \emptyset  \quad \forall 1\leq i < j <k \leq \ell, \; \forall e \in E.
\end{equation}
\end{definition}


The key differences between definition of a simple cycle and an $\alpha$-cycle are that condition \eqref{eq alpha diff} is omitted, and the indices in condition \eqref{eq alpha cond} are no longer required to be contiguous.
While these differences may appear minor, they play a crucial role in our proof of \cref{th main}.
We are now ready to state our new characterization for $\alpha$-acyclic hypergraphs.

\begin{theorem}
\label{th characterization}
A hypergraph $G$ is $\alpha$-acyclic if and only if it does not contain a simple cycle.
\end{theorem}

To prove~\cref{th characterization}, we make use of the following lemmata relating simple cycles to those $\alpha$-cycles that cannot be ``shortened.''
We first show that a simple cycle is an $\alpha$-cycle.

\begin{lemma}
\label{lem relate}
A simple cycle is an $\alpha$-cycle.
\end{lemma}
\begin{prf}
    Let $G=(V,E)$ be a hypergraph and let $C=e_1, e_2, \ldots, e_{\ell}, e_{\ell+1}$ with $e_{\ell+1} = e_1$ and for some $\ell \geq 3$ be a simple cycle of $G$. We show that $C$ satisfies conditions~\eqref{eq alpha diff} and~\eqref{eq alpha cond} implying that it is an $\alpha$-cycle. We first consider~\eqref{eq alpha cond}. Letting $i=j-1$ and $k=j+1$ in~\eqref{eq simple} we deduce that $C$ satisfies
    condition~\eqref{eq alpha cond} for $2 \leq i \leq \ell-1$.
    Letting $i=1$, $j=2$ and $k=\ell$ (resp. $i=1$, $j=\ell-1$ and $k=\ell$) in~\eqref{eq simple} we deduce that $C$ satisfies
    condition~\eqref{eq alpha cond} for $i=1$ (resp. $i=\ell$). Therefore, $C$ satisfies~\eqref{eq alpha cond}.
    
    Next, consider condition~\eqref{eq alpha diff}; that is, for any $i\neq j \in [\ell]$ we want to show that $e_i\setminus e_j \neq \emptyset$.
    Recall that for any $i \in [\ell]$, we have $e_i \supseteq s_{i-1} \cup s_i$. Let $q = \min(i,j)$ and $p=\max(i,j)$.
%
First assume that $q >1$; in this case letting $i=q-1$, $j=q$, $k=p$ and $e=e_q$ in
        ~\eqref{eq simple}, we obtain $s_p \setminus e_q \neq \emptyset$, which implies that $e_p \setminus e_q \neq \emptyset$ for all $2\leq q < p \leq \ell$. Next, let $q=1$; in this case letting $i=1$, $j=p$, $k=\ell$
        and $e=e_1$ in~\eqref{eq simple} we obtain $s_p \setminus e_1 \neq \emptyset$, which implies that $e_p \setminus e_1 \neq \emptyset$ and $e_{p+1} \setminus e_1 \neq \emptyset$ for all $1 < p \leq \ell-1$. Therefore, 
        $e_p \setminus e_1 \neq \emptyset$ for all $p \in [\ell] \setminus \{1\}$.
Hence, $C$ satisfies condition~\eqref{eq alpha diff} and as a result it is a simple cycle.
\end{prf}

We say that an $\alpha$-cycle $C= e_1, e_2, \ldots, e_{\ell}, e_{\ell+1}$ is \emph{chordless} if for every $1 \leq i < j \leq \ell$ with $j-i \le \ell-3$ and for every $\tilde e \in E \setminus \{e_i,\ldots e_j\}$, the edge sequence $e_i, e_{i+1}, \ldots, e_j, \tilde e, e_i$ is not an $\alpha$-cycle. 
Clearly, an $\alpha$-cycle of length three is chordless. As for chordless cycles in graphs, chordless $\alpha$-cycles are  $\alpha$-cycles that cannot be ``shortened.''
We next prove a property of chordless $\alpha$-cycles which relates them to simple cycles:

\begin{lemma}\label{diagonals}
Let $G=(V,E)$ be a hypergraph and let $C=e_1,e_2,\ldots,e_{\ell},e_{\ell+1}$, for some $\ell \geq 4$, be an $\alpha$-cycle of $G.$ Define $s_i = e_i \cap e_{i+1}$ for $i \in [\ell]$.
If $C$ is chordless, then for every $e \in E$, we have
\begin{align}
\tag{C}
\label{eq chordless}
(s_i \cup s_j) \setminus e \neq \emptyset & \quad \forall i<j \in [\ell]: j-i > 1, \; (i,j)\neq (1,\ell), \; \forall e \in E.
\end{align}
\end{lemma}
\begin{prf}
Denote by $\tilde e$ a maximal edge of $G$ that does not satisfy~\eqref{eq chordless}. Notice that if there is an edge in $E$ that violates~\eqref{eq chordless}, then any maximal edge of $G$ containing this edge also violates~\eqref{eq chordless}. Denote by $L$ the set of indices $i \in [\ell]$ for which we have $\tilde e \supseteq s_i$.
We then have $\tilde e \supseteq \cup_{i \in L}{s_i}$.
Notice that
we have $|L| \geq 2$.
Let $i_1$ be the smallest index in $L$ and let $i_2$ be the second smallest.
If $i_2 - i_1 > 1$, then we set $p=i_1$, $q=i_2$. Notice that by~\eqref{eq chordless}, we must have $(i_1, i_2) \neq (1, \ell)$.
If $i_2 - i_1 = 1$, then $L$ should contain at least three elements because by assumption $\tilde e$ violates assumption~\eqref{eq chordless}.
Denote by $i_3$ the third smallest index in $L$.  We then set $p=i_2$, $q=i_3$.
In this case, since $C$ is an $\alpha$-cycle, by~\eqref{eq alpha cond} we must have $i_3-i_2 > 1$.

%
We claim that the edge sequence $C'=\tilde e, e_{p+1}, e_{p+2}, \ldots, e_q, \tilde e$ is an $\alpha$-cycle, which contradicts the assumption that $C$ is chordless.
Since $\tilde e$ is a maximal edge of $G$, to verify condition~\eqref{eq alpha diff} it suffices to show that $e_i \setminus \tilde e \neq \emptyset$ for all $i \in \{p+1, \ldots, q\}$.
This follows since, by definition of $L$, $\tilde e$ does not contain any $s_i$ for $p < i < q$ and we have $s_i = e_i \cap e_{i+1}$.
Notice that since by assumption $q-p > 1$, $C'$ has at least three distinct edges.
Then, to prove that $C'$ is an $\alpha$-cycle, it suffices to check if condition~\eqref{eq alpha cond} holds.
We consider two cases:
\begin{itemize}
\item If $q - p > 2$, then we need to show that for every $e \in E$:
$$
((\tilde e \cap e_{p+1}) \cup s_{p+1} \cup s_{p+2})\setminus e \neq \emptyset, \quad (s_{q-2}\cup s_{q-1}\cup(e_q \cap \tilde e)) \setminus e \neq \emptyset.
$$
Since by assumption  $\tilde e \supseteq s_p \cup s_q$, it follows that $\tilde e \cap e_{p+1} \supseteq s_p$, and $e_q \cap \tilde e \supseteq s_q$. Since $C$ is an $\alpha$-cycle, by~\eqref{eq alpha cond}, the above relations are valid.

\item If $q-p = 2$, then we need to show that for every $e\in E$:
$$
((\tilde e \cap e_{p+1}) \cup s_{p+1} \cup (e_q \cap \tilde e) )\setminus e \neq \emptyset,
$$
whose validity follows from~\eqref{eq alpha cond} since
$\tilde e \cap e_{p+1} \supseteq s_p$, $e_q \cap \tilde e \supseteq s_q$, and $C$ is an $\alpha$-cycle.
\end{itemize}
\end{prf}

We then obtain the following.

\begin{lemma}\label{chordSimple}
A chordless $\alpha$-cycle is a simple cycle.
\end{lemma}

\begin{prf}
Let $G=(V,E)$ be a hypergraph and let $C=e_1,e_2,\ldots,e_{\ell}, e_{\ell+1}$, with $e_{\ell+1} = e_1$, be a chordless $\alpha$-cycle of $G$. 
If $\ell = 3$, the definitions of $\alpha$-cycles, chordless $\alpha$-cycles, and simple cycles coincide. Hence, let $\ell \geq 4$. We show that any chordless $\alpha$-cycle $C$ satisfies condition~\eqref{eq simple}, implying that it is a simple cycle.
 Consider any $i,j,k \in [\ell]$ satisfying $i < j < k$.
Notice that $k-i > 1$.
If $(i,k) \neq (1, \ell)$, then the statement follows since by~\eqref{eq chordless} we have $(s_i \cup s_k) \setminus e \neq \emptyset$ for every $e \in E$.
Now suppose that $(i,k) = (1, \ell)$. Since $i < j < k$ and $\ell \geq 4$, we either have $j-i \geq 2$ or $k-j \geq 2$, in which case by~\eqref{eq chordless} we get $(s_i \cup s_j) \setminus e \neq \emptyset$ for all $e \in E$ or $(s_j \cup s_k)\setminus e \neq \emptyset$, respectively. Hence, condition~\eqref{eq simple} holds for $C$.
\end{prf}

\begin{figure}
\centering
\includegraphics[scale=0.6]{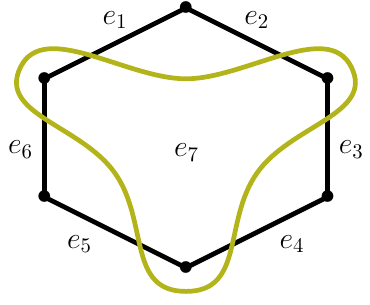}
\caption{The edge sequence $C=e_1, e_2, \cdots, e_6,e_1$ is an $\alpha$-cycle. However, $C$ is not a simple cycle because $e_7 \supseteq s_2 \cup s_4 \cup s_6$. The sequence $e_1,e_2,e_7,e_1$ is a simple cycle.}
\label{fig2}
\end{figure}

\begin{figure}
\centering
\includegraphics[scale=0.6]{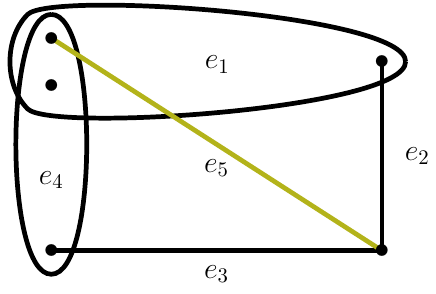}
\caption{The edge sequence $C=e_1, e_2, e_3, e_4, e_1$ is a simple cycle which also satisfies condition~\eqref{eq chordless}. However, $C$ is not a chordless $\alpha$-cycle because it can be shortened; \ie $e_1, e_2, e_5, e_1$ is a shorter $\alpha$-cycle.}
\label{fig3}
\end{figure}

Figures~\ref{fig2} and~\ref{fig3} illustrate the difference between $\alpha$-cycles, chordless $\alpha$-cycles, and simple cycles. Thanks to~\cref{lem relate} and~\cref{chordSimple}, we can now easily prove \cref{th characterization}.

\paragraph{Proof of~\cref{th characterization}.} 
If a hypergraph contains an $\alpha$-cycle, then it contains a chordless $\alpha$-cycle and by~\cref{chordSimple}, it contains a simple cycle. If a hypergraph does not contain an $\alpha$-cycle,  by~\cref{lem relate}, it does not contain a simple cycle either.
\qed

\section{Sufficiency of $\alpha$-acyclicity}
\label{sec sufficient}

In this section, we prove the ``easy'' direction of \cref{th main}:

\begin{theorem}
\label{th sufficiency}
Let $G=(V,E)$ be an $\alpha$-acyclic hypergraph.
Then $\CER(G)$ is an extension of $\MP(G)$.
\end{theorem}

As we discussed earlier, \cref{th sufficiency} is implied by the results in~\cite{Lau09,BieMun18} and \cite{dPKha21MOR}.
In~\cite{Lau09,BieMun18}, the authors give extended formulations for the convex hull of the feasible set of (possibly constrained) binary polynomial optimization problems. The size of these extended formulations is parameterized in terms of the ``treewidth'' of their so-called intersection graphs. As we detail in~\cite{dPKha21MOR}, in the unconstrained case, their result can be equivalently stated as: Let $G$ be an $\alpha$-acyclic hypergraph of rank $r$. Then $\MP(G)$ has an extended formulation with at most $2^r|E_{\max}|$ variables and inequalities.  While not stated explicitly, the proofs in~\cite{Lau09,BieMun18} imply that $\CER(G)$ is an extended formulation for $\MP(G)$.

Here, we present a simple proof that does not involve treewidth and makes direct use of $\alpha$-acyclicity.
Our proof relies on a result of~\cite{dPKha18MPA} regarding the decomposability of the multilinear polytope.
To state this result, we introduce some terminology.
Given a hypergraph $G=(V,E)$, and $V' \subseteq V$, the \emph{section hypergraph} of $G$ induced by $V'$ is the hypergraph $G'=(V',E')$, where $E' = \{e \in E : e \subseteq V'\}$.
Given hypergraphs $G_1=(V_1,E_1)$ and $G_2=(V_2,E_2)$, we denote by $G_1 \cap G_2$ the hypergraph $(V_1 \cap V_2, E_1 \cap E_2)$,
and we denote by $G_1 \cup G_2$, the hypergraph $(V_1 \cup V_2, E_1 \cup E_2)$. 
In the following, we consider a hypergraph $G$, and two distinct section hypergraphs of $G$, denoted by $G_1$ and $G_2$, such that $G_1 \cup G_2 = G$.
We say that $\MP(G)$ is \emph{decomposable} into $\MP(G_1)$ and $\MP(G_2)$ if 
the system comprised of a description of $\MP(G_1)$ and a description of $\MP(G_2)$, is a description of $\MP(G)$.

\begin{theorem}[theorem~1 in~\cite{dPKha18MPA}]
\label{th decomposition}
Let $G_1$,$G_2$ be section hypergraphs of a hypergraph $G$ such that $G_1 \cup G_2 = G$ and $G_1 \cap G_2$ is a complete hypergraph.
Then, $\MP(G)$ is decomposable into $\MP(G_1)$ and $\MP(G_2)$.
\end{theorem}



\begin{figure}
    \centering
    \hspace{-1.2cm}
    \hspace{\stretch{1}}
    \includegraphics[scale=.55]{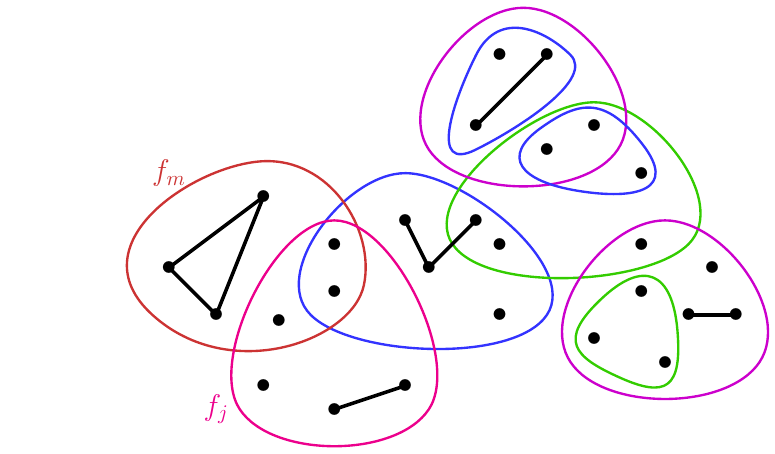}
    \hspace{\stretch{3}}
    \includegraphics[scale=.55]{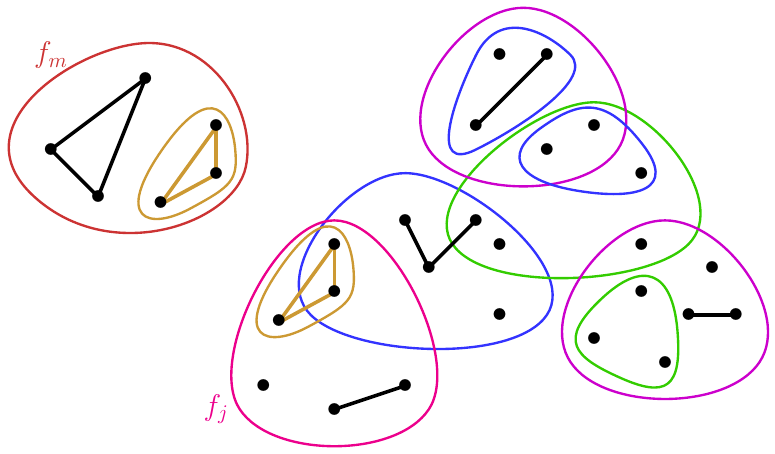}
    \hspace{\stretch{1}}
    \caption{An illustration of the proof technique for~\cref{th sufficiency}.}
    \label{fig alpha}
\end{figure}

In the following, we will often use the fact that, for hypergraphs, adding or removing edges inside a maximal edge preserves $\alpha$-acyclicity.
This follows, for example, from the definition of acyclicity (page~483) and theorem~3.4 in \cite{BeeFagMaiYan83}.
In particular, this fact implies that the hypergraph $G$ is $\alpha$-acyclic if and only if its completion $\cl(G)$ is $\alpha$-acyclic.
We are now ready to present our proof. See Figure~\ref{fig alpha} for an illustration of our proof technique.

\paragraph{Proof of \cref{th sufficiency}.}
To prove \cref{th sufficiency} we show that, if the hypergraph $G=(V,E)$ satisfying $G = \cl(G)$ is $\alpha$-acyclic, then $\CER(G)=\MP(G)$.
We prove this statement by induction on $|\ME|$, where $\ME$ denotes the set of maximal edges of $G$.
The base $|\ME|=1$ case follows from \cref{prop:RLT}.

In the remainder of the proof we show the step case, and we assume $|\ME| \ge 2$.
The hypergraph $(V,\ME)$ is $\alpha$-acyclic.
By definition of $\alpha$-acyclicity, there is a running intersection ordering of $\ME$, say $f_1, f_2, \ldots, f_m$.
In particular, there exists $j < m$, such that $f_m \cap\left(f_1 \cup f_2 \cup \ldots \cup f_{m-1}\right) \subseteq f_j$.
Let $G_1$ be the section hypergraph of $G$ induced by $f_m$.
The hypergraph $G_1$ is complete, therefore from \cref{prop:RLT} implies $\CER(G_1)=\MP(G_1)$.

Let $G_2$ be the section hypergraph of $G$ induced by $f_1 \cup f_2 \cup \ldots \cup f_{m-1}$.
We claim that the set of maximal edges of $G_2$ is $\ME \setminus \bra{f_m}$.
Clearly, $f_1, f_2, \ldots, f_{m-1}$ are maximal edges of $G_2$.
Since every edge of $G_2$ is also an edge of $G$, we only need to show that no edge $f$ of $G_2$ with $f \subseteq f_m$ is a maximal edge of $G_2$.
So let $f$ be an edge of $G_2$ with $f \subseteq f_m$.
Then $f \subseteq f_1 \cup f_2 \cup \ldots \cup f_{m-1}$.
Hence, $f = f \cap\left(f_1 \cup f_2 \cup \ldots \cup f_{m-1}\right) \subseteq f_m \cap\left(f_1 \cup f_2 \cup \ldots \cup f_{m-1}\right) = f_j$.
Since $f_j$ is an edge of $G_2$, $f$ is not a maximal edge of $G_2$.

We also observe that $G_2$ is $\alpha$-acyclic, since the running intersection ordering $f_1, f_2, \ldots, f_m$ of $\ME$ directly implies that $f_1, f_2, \ldots, f_{m-1}$ is a running intersection ordering of $\ME \setminus \bra{f_m}$, which we showed is the set of maximal edges of $G_2$.
So $G_2$ is $\alpha$-acyclic, has $m-1$ maximal edges, and clearly $G_2 = \cl(G_2)$.
Thus by induction $\CER(G_2) = \MP(G_2)$.

Observe that $G_1 \cup G_2 = G$, since every node of $G$ is contained in at least one edge of $G$, and each edge of $G$ is contained in at least one edge among $f_1, f_2, \dots, f_m$.
Furthermore, $G_1 \cap G_2$ is a complete hypergraph, due to the definition of $G$, since $G_1$ is complete.
We can now apply \cref{th decomposition} and obtain that $\MP(G)$ is decomposable into $\MP(G_1)$ and $\MP(G_2)$.
Since $\CER(G_1) = \MP(G_1)$ and $\CER(G_2) = \MP(G_2)$, the system comprised of a description of $\CER(G_1)$ and a description of $\CER(G_2)$, that is precisely the description of $\CER(G)$, is a description of $\MP(G)$. \qed
\medskip

In~\cite{dPKha18SIOPT}, we give a polynomial-size extended formulation of the multilinear polytope of $\gamma$-acyclic hypergraphs with at most $|V| + 2|E|$ variables and at most $|V| + (r +2)|E|$ inequalities. In~\cite{dPKha23mMPA}, we further generalize this result by obtaining a polynomial-size extended formulation of the multilinear polytope of $\beta$-acyclic hypergraphs with at most $(r-1) |V|+|E|$ variables and at most
$(3r-4)|V|+4|E|$ inequalities. 
By~\cref{th sufficiency}, the complete edge relaxation is also an extended formulation of the multilinear polytope of $\gamma$-acyclic and $\beta$-acyclic hypergraphs, 
because every $\gamma$-acyclic or $\beta$-acyclic hypergraph is $\alpha$-acyclic. However, the extended formulations of~\cite{dPKha18SIOPT,dPKha23mMPA} grow linearly in the rank $r$, while the size of the complete edge relaxation is exponential in the rank. Recall that $\CER(G)$ is defined by at most $2^r|E_{\max}|$ variables and inequalities.
Therefore, for the multilinear polytope of $\gamma$-acyclic or $\beta$-acyclic hypergraphs with unbounded rank, the aforementioned extended formulations are of polynomial-size while $\CER(G)$ is not.
It is worth mentioning that $\CER(G)$ is simpler to construct in the sense that it requires only the knowledge of the set of maximal edges of $G$. 
On the other hand, the extended formulation in~\cite{dPKha18SIOPT} requires $G$ to be $\gamma$-acyclic as it uses laminar hypergraphs as its building block, and the extended formulation in~\cite{dPKha23mMPA} requires $G$ to be $\beta$-acyclic together with knowledge of a nest-point elimination ordering.

\section{Necessity of $\alpha$-acyclicity}
\label{sec necessary}

In this section, we prove the main result of this paper:

\begin{theorem}
\label{th necessity}
Let $G=(V,E)$ be a hypergraph and suppose that $\CER(G)$ is an extension of $\MP(G)$. Then $G$ is an $\alpha$-acyclic hypergraph.
\end{theorem}

To prove~\cref{th necessity}, we show that if $G$ has a simple cycle, then $\CER(G)$ is not an extension of $\MP(G)$. To this end,  we first introduce some operations on hypergraphs; namely, node fixing, node contraction, and node expansion. We show how the multilinear polytope and the complete edge relaxation are transformed under these operations. Through a recursive application of these operations, we then prove that one can reduce the length of simple cycles. This in turn implies that we can limit our attention to hypergraphs containing simple cycles of length three, a structure that can be tackled directly.

\subsection{The node fixing operator}

Consider a hypergraph $G = (V,E)$ and let $V' \subseteq V$. The \emph{subhypergraph} of $G$
induced by $V'$ is the hypergraph $G_{V'}$  with node set 
$V'$ and with edge set $E_{V'}:=\{e \cap V' : e \in E, \; |e\cap V'| \geq 2\}$.
For every edge $e$ of $G_{V'}$, there may exist several edges $e'$ of $G$ satisfying $e = e'\cap V$; we denote by
$e'(e)$ one such arbitrary edge of $G$. For ease of notation, we identify an edge $e$ of $G_{V'}$ with an
edge $e'(e)$ of $G$. Define
\begin{equation}\label{affineL}
L_{V'} := \Big\{z \in \R^{V\cup E} : z_v = 1 \ \forall v \in V \setminus V'\Big\}.
\end{equation}
Let $R(G)$ denote a convex relaxation of $\MP(G)$.
Denote by $\proj_{G_{V'}} (R(G) \cap L_{V'})$ the set obtained from $R(G) \cap L_{V'}$ by projecting out all variables $z_v$, for $v \in V \setminus V'$, and
$z_{e}$, for $e \in E \setminus E'$, where $E':=\{e'(e) \in E: e \in E_{ V'}\}$.
Define 
\begin{equation}\label{affinecL}
\bar L_{V'} := \Big\{z \in \R^{V\cup\cl(E)} : z_v = 1 \ \forall v \in V \setminus V'\Big\}.
\end{equation}
Let $R(\cl(G))$ denote a convex relaxation of $\MP(\cl(G))$.
Denote by $\proj_{G_{V'}} (R(\cl(G)) \cap \bar L_{V'})$ the set obtained from $R(\cl(G)) \cap \bar L_{V'}$ by projecting out all variables $z_v$ for all $v \in V \setminus V'$, and $z_f$ for all $f \subseteq e$ such that $f \not\subseteq V'$ all $e \in E$.

The next lemma indicates that the multilinear polytope (resp. the complete edge relaxation) of a subhypergraph, \ie $\MP(G_{V'})$ (resp. $\CER(G_{V'})$) can be described as a projection of the face of $\MP(G)$ (resp. $\CER(G)$) defined by the affine set $L_{V'}$ (resp. $\bar L_{V'}$).

\begin{lemma}
\label{res0}
Let $G = (V, E)$ be a hypergraph and let $L_{V'}$ and $\bar L_{V'}$ with $V'=V \setminus \{v'\}$ for some $v' \in V$ be defined by~\eqref{affineL} and~\eqref{affinecL}, respectively. Then we have:
\medskip
\begin{itemize}
    \item [(i)] $\MP(G_{V'})= \proj_{G_{V'}} (\MP(G) \cap L_{V'})$.
    \item [(ii)] $\CER(G_{V'}) = \proj_{G_{V'}} (\CER(G) \cap \bar L_{V'})$.
\end{itemize}
\end{lemma}

\begin{prf}
Let $\S(G)$ denote the multilinear set associated with $G$.
First, consider Part~(i). Denote by $\proj_{G_{V'}} (\S(G) \cap L_{V'})$  the set obtained from $\S(G) \cap L_{V'}$ by projecting out variables $z_{v'}$, and $z_{e}$, $e \in E \setminus E'$.
From the definition of the subhypergraph it follows that $\S(G_{V'}) = \proj_{G_{V'}} (\S(G) \cap L_{V'})$, implying that $\MP(G_{V'})= \proj_{G_{V'}} \conv(\S(G) \cap L_{V'})$. Moreover, the set $\MP(G) \cap L_{V'}$ is a face of $\MP(G)$ defined by the valid inequality $z_{v'} \leq 1$, implying that $\MP(G) \cap L_{V'} = \conv(\S(G)\cap L_{V'})$. We then deduce that $\MP(G_{V'}) = \proj_{G_{V'}}(\MP(G) \cap L_{V'})$.

\medskip

Next, consider Part~(ii). We start by characterizing the polytope $\proj_{G_{V'}} (\CER(G) \cap \bar L_{V'})$ and then show that it coincides with $\CER(G_{V'})$. Denote by $\ME$ the set of maximal edges of $G$. 
Recall that $\CER(G)$ is obtained by putting together the inequalities defining $\MP(K_{\bar e})$ for all $\bar e \in \ME$, where $K_{\bar e}$ denotes the complete hypergraph with node set $\bar e$ and $\MP(K_{\bar e})$ is obtained using~\cref{prop:RLT}: 
$$\psi(U, \bar e) \geq 0, \quad \forall U \subseteq \bar e,$$
and where $\psi(U, \bar e)$ is defined by~\eqref{defRLT}. Now consider any edge $\bar e \in \ME$ such that $\bar e \ni v'$, and let $\bar U$ denote the set of subsets of $\bar e$ such that for any $U \in \bar U$, we have $U \ni v'$. Since by definition $z_e = \prod_{v \in e}{z_v}$ for all $e \subseteq \bar e$, it follows that  $\psi(U, \bar e) = \prod_{v \in U} {(1-z_v)} \prod_{v \in \bar e \setminus U}{z_v}$. 
We then deduce that
$$\sum_{U \in \bar U} {\psi(U, \bar e)} = 1-z_{v'} = 0,$$ 
which together with $\psi(U,\bar e) \geq 0$ for all $U \in \bar U$ implies that
\begin{equation}\label{simplified}
\psi(U,\bar e) = 0, \quad \forall U \in \bar U.
\end{equation}
First consider $U = \{v'\}$. In this case~\eqref{simplified} simplifies to: 
\begin{equation}\label{simple1}
z_{\bar e} = z_{\bar e \setminus \{v'\}}.
\end{equation}
Then for any $v \in \bar e \setminus \{v'\}$, consider $U = \{v, v'\}$. Substituting in~\eqref{simplified} and using~\eqref{simple1} we obtain:
\begin{equation}\label{simple2}
z_{\bar e \setminus \{v\}}= z_{\bar e \setminus \{v,v'\}}, \quad \forall v \in \bar e \setminus \{v'\}.
\end{equation}
Next for any $u \neq v \in \bar e \setminus \{v'\}$, consider $U = \{u, v, v'\}$. Substituting in~\eqref{simplified} and using~\eqref{simple1} and~\eqref{simple2} yields $z_{\bar e \setminus \{u, v\}}= z_{\bar e \setminus \{u, v,v'\}}$ for all $u \neq v \in \bar e \setminus \{v'\}$. Applying this argument recursively, we deduce that for each $\bar e \in \ME$ with $\bar e \ni v'$ the polytope $\MP(K_{\bar e}) \cap \{z: z_{v'=1}\}$ is given by:
\begin{align}
    & z_{e} = z_{e\setminus \{v'\}}, \quad \forall e \subseteq \bar e: \; |e| \geq 2, \; e\ni v' \label{simpf1}\\
    & \psi(U, \bar e \setminus \{v'\}) \geq 0, \quad \forall U \subseteq \bar e \setminus \{v'\}\label{simpf2}.
\end{align}
Notice that by~\cref{prop:RLT}, inequalities~\eqref{simpf2} define the polytope $\MP(K_{\bar e \setminus \{v'\}})$, where $K_{\bar e \setminus \{v'\}}$ denotes the complete hypergraph with node set $\bar e \setminus \{v'\}$.
Moreover, the relaxation $\CER(G_{V'})$ is obtained by putting together the descriptions of $\MP(K_{\bar e \setminus \{v'\}})$ for all edges $\bar e \in \ME$ with $\bar e \ni v'$ and
$\MP(K_{\bar e})$ for all edges $\bar e \in \ME$ with $\bar e \not\ni v'$.
Hence, using equalities~\eqref{simpf1} to project out variables $z_e$, for any edge $e \ni v'$, we obtain $\CER(G_{V'}) = \proj_{G_{V'}} (\CER(G) \cap \bar L_{V'})$.

\end{prf}



Thanks to~\cref{res0}, we can show that the tightness of the complete edge relaxation is preserved under the node fixing operation.

\begin{proposition}\label{cor: fixing}
    Let $G=(V,E)$ be a hypergraph and let $\bar V \subseteq V$.
    If $\CER(G)$ is an extension of $ \MP(G)$, then
    $\CER(G_{\bar V})$ is an extension of $\MP(G_{\bar V})$.
\end{proposition}

\begin{prf}
Consider the affine set $L_{V'}$ (resp. $\bar L_{V'}$) defined by~\eqref{affineL} (resp.~\eqref{affinecL} ) with $V' = V \setminus \{v'\}$ for some $v' \in V$. Denote by $\proj_G(\CER(G))$ the set obtained from $\CER(G)$ by projecting out variables $z_e$, for all $e \in \cl(E) \setminus E$.
Since $\proj_G(\CER(G))=\MP(G)$, we have $\proj_G(\CER(G)) \cap L_{V'}= \proj_G(\CER(G) \cap \bar L_{V'}) =\MP(G) \cap L_{V'}$, where the first equality follows since the operator $\proj_{G}(\cdot)$ does not project out variables $z_v$, $v \in V \setminus V'$. 
This in turn implies that 
$$\proj_{G_{V'}}(\proj_G(\CER(G) \cap \bar L_{V'})) = \proj_{G}(\proj_{G_{V'}}(\CER(G) \cap \bar L_{V'})) = \proj_{G_{V'}}(\MP(G) \cap L_{V'}).$$
Therefore, by Parts~(i) and~(ii) of~\cref{res0}, we deduce that $\proj_{G}(\CER(G_{V'})) = \MP(G_{V'})$. The proof then follows from a recursive application of this argument. 
\end{prf}

\subsection{The node contraction operator}

Let $G=(V,E)$ be a hypergraph, and let $u,w \in V$ such that there exists $e \in E$ with $e \ni u,w$. 
We say that the hypergraph $G_{w\rightarrow u}=(V_{w\rightarrow u}, E_{w\rightarrow u})$ is obtained from $G$ by \emph{contracting} $w$ to $u$, if $V_{w\rightarrow u} := V \setminus \bra{w}$, and
$$
E_{w\rightarrow u} := \bra{e \in E: e \not\ni w} \cup \bra{(e \setminus \bra{w}) \cup \bra{u} : e \in E \setminus \{\{u,w\}\}, \ e \ni w}.
$$
Define the subspace:
\begin{equation}\label{newaff}
    Q_{w \rightarrow u} := \Big\{z \in \R^{V \cup E}:  z_w = z_u \Big\}.
\end{equation}
Let $R(G)$ denote a convex relaxation of $\MP(G)$.
 Denote by $\proj_{G_{w \rightarrow u}} (R(G) \cap Q_{w \rightarrow u})$   the set obtained from $R(G) \cap Q_{w\rightarrow u}$ by projecting out variables $z_w$ and $z_e$, $e \in E$ such that $e=\{u,w\}$ or $e \ni w$ and $(e \setminus \{w\}) \cup \{u\} \in E$.
 Define
 \begin{equation}\label{newaffc}
    \bar Q_{w \rightarrow u} := \Big\{z \in \R^{V \cup \cl(E)}:  z_w = z_u \Big\}.
\end{equation}
Let $R(\cl(G))$ denote a convex relaxation of $\MP(\cl(G))$.
Denote by $\proj_{G_{w \rightarrow u}} (R(\cl(G)) \cap \bar Q_{w \rightarrow u})$ the set obtained from $R(\cl(G)) \cap \bar Q_{w\rightarrow u}$ 
by projecting out variables $z_w$ and $z_e$ for all $e \in \cl(E)$ such that $e=\{u,w\}$ or $e \ni w$ and $(e \setminus \{w\}) \cup \{u\} \in \cl(E)$. 
We then have the following result.

The next lemma indicates that the multilinear polytope of the hypergraph obtained by contracting node $w$ to $u$, \ie  $\MP(G_{w\rightarrow u})$  can be described as a projection of the face of $\MP(G)$ defined by $z_w = z_u$. Moreover, the complete edge relaxation $\CER(G_{w\rightarrow u})$ is contained in a projection of the face of $\CER(G)$ defined by $z_w = z_u$.

\begin{lemma}\label{lemma:contract}
Let $G = (V, E)$ be a hypergraph, let $u,w \in V$, let $G_{w\rightarrow u}$ be the hypergraph obtained from $G$ by contracting $w$ to $u$, and let $Q_{w\rightarrow u}$ and $\bar Q_{w\rightarrow u}$ be defined by~\eqref{newaff} and~~\eqref{newaffc}, respectively. We then have:
\medskip
\begin{itemize}
    \item [(i)] $\MP(G_{w\rightarrow u})= \proj_{G_{w\rightarrow u}} (\MP(G) \cap Q_{w\rightarrow u})$.
    \item [(ii)] $\CER(G_{w\rightarrow u}) \subseteq \proj_{G_{w\rightarrow u}} (\CER(G) \cap \bar Q_{w\rightarrow u})$.
\end{itemize}
\end{lemma}

\begin{prf}
Let $\S(G)$ denote the multilinear set associated with $G$.
First, consider Part~(i). Define $\tilde G := (V, \tilde E)$, where $\tilde E := E \cup\{\{u,w\}\}$. Clearly, if $\{u,w\} \in E$, then the two hypergraphs $G$ and $\tilde G$ coincide. 
Define $\tilde Q_{w\rightarrow u} := \{z\in \R^{V \cup \tilde E}: z_w = z_u\}$.
Denote by $\proj_{\tilde G_{w\rightarrow u}} (\S(\tilde G) \cap \tilde Q_{w\rightarrow u})$  the set obtained from $\S(\tilde G) \cap \tilde Q_{w\rightarrow u}$ by projecting out variables $z_w$ and $z_e$, $e \in E$ such that $e=\{u,w\}$ or $e \ni w$ and $(e\setminus \{w\})\cup \{u\} \in E$. 
From the definition of the contraction operator it follows that $\S(G_{w\rightarrow u}) = \proj_{\tilde G_{w\rightarrow u}} (\S(\tilde G) \cap \tilde Q_{w\rightarrow u})$, implying that $\MP(G_{w\rightarrow u})= \proj_{\tilde G_{w\rightarrow u}} \conv(\S(\tilde G) \cap \tilde Q_{w\rightarrow u})$. Moreover, the inequality $z_u + z_w \geq 2 z_{\{u,w\}}$
is valid for $\S(\tilde G)$ and is satisfied tightly at a binary point $z$ only if $z_w = z_u$. It then follows that $\MP(\tilde G) \cap \tilde Q_{w\rightarrow u}$ is a face of $\MP(\tilde G)$ defined by this inequality, and therefore we have $\MP(\tilde G) \cap \tilde Q_{w\rightarrow u} = \conv(\tilde \S(G)\cap \tilde Q_{w\rightarrow u})$. We then deduce that 
\begin{equation}\label{auxfp}
\MP(G_{w\rightarrow u}) = \proj_{\tilde G_{w\rightarrow u}}(\MP(\tilde G) \cap \tilde Q_{w\rightarrow u})
\end{equation} 
Now consider $\MP(G)$; if $\{u,w\} \in E$, then $\MP(\tilde G) = \MP(G)$ and we are done. Otherwise, $\MP(G)$ is obtained by projecting out the variable $z_{\{u,w\}}$ from the description of $\MP(\tilde G)$. Let us denote this operation as $\MP(G) = \proj_{G}(\MP(\tilde G))$. We then have 
\begin{align*}
\proj_{ G_{w\rightarrow u}}(\MP(G) \cap Q_{w\rightarrow u}) = & \proj_{G_{w\rightarrow u}}(\proj_G(\MP(\tilde G))\cap Q_{w\rightarrow u}) \\
= & \proj_{G_{w\rightarrow u}}(\proj_G(\MP(\tilde G)\cap \tilde Q_{w\rightarrow u}))\\
= & \proj_{\tilde G_{w\rightarrow u}}(\MP(\tilde G)\cap \tilde Q_{w\rightarrow u}), 
\end{align*}
where the second equality follows since the operator $\proj_G (\cdot)$ does not project out variables $z_u, z_w$. Together with~\eqref{auxfp}, this in turn implies that
$\MP(G_{w\rightarrow u}) = \proj_{G_{w\rightarrow u}}(\MP(G) \cap Q_{w\rightarrow u})$.

\medskip

Next, consider Part~(ii). We first define a polyhedral relaxation of $\S(G) \cap Q_{w\rightarrow u}$ in the extended space $\R^{V\cup\cl(E)}$, which we denote by $\overline{\CER}(G)$, satisfying $\overline{\CER}(G) \subseteq \CER(G)\cap \bar Q_{w\rightarrow u}$. Subsequently, we show that $\CER(G_{w\rightarrow u}) = \proj_{G_{w\rightarrow u}}(\overline{\CER}(G))$. This in turn implies that $\CER(G_{w\rightarrow u}) \subseteq \proj_{G_{w\rightarrow u}} (\CER(G) \cap \bar Q_{w\rightarrow u})$.

Let us start by defining $\overline{\CER}(G)$; Denote by $\ME$ the set of maximal edges of $G$. For each $\bar e \in \ME$, as before denote by $K_{\bar e}$ the complete hypergraph with the node set $\bar e$. Then $\overline{\CER}(G)$ is obtained by intersecting the following polyhedral descriptions for all $\bar e \in \ME$:

\medskip
\begin{enumerate}
    \item For each  $\bar e \in \ME$ such that $\bar e \not\ni w$, the multilinear polytope $\MP(K_{\bar e})$ obtained using~\cref{prop:RLT}.

    \item For each $\bar e \in \ME$ such that $\bar e \ni w$ and $\bar e \not\ni u$, the multilinear polytope $\MP(K_{\bar e})$ obtained using~\cref{prop:RLT} together with the equality constraint $z_w = z_u$.

    \item For each $\bar e \in \ME$ such that $\bar e \supseteq \{u,w\}$, the convex hull of the set 
    
    $$\S^{=}(K_{\bar e}): =\Big\{z \in \{0,1\}^{P(\bar e) \setminus \{\emptyset\}}: z_e = \prod_{v\in e}{z_v}, \; \forall e \subseteq \bar e, \; |e|\geq 2, \; z_u = z_w\Big\},$$ 
    
    using RLT~\cite{SheAda90}, where $P(\bar e)$ denotes the power set of $\bar e$.

\end{enumerate}
\medskip

Notice that the difference between $\CER(G) \cap \bar Q_{w\rightarrow u}$ and $\overline{\CER}(G)$ is in Case~(3) above.
That is, while to construct $\CER(G) \cap \bar Q_{w\rightarrow u}$ for each maximal edge $\bar e \supseteq \{u,w\}$, we first construct $\MP(K_{\bar e})$ using RLT and then intersect this polyhedron with $z_u = z_w$, to construct $\overline{\CER}(G)$, we directly construct the convex hull of the set $\S^{=}(K_{\bar e})$, which can be written as $\S^{=}(K_{\bar e}) = \S(K_{\bar e})\cap \{z\in \{0,1\}^{P(\bar e) \setminus \{\emptyset\}}: z_u = z_w\}$, using RLT. It then follows that $\overline{\CER}(G) \subseteq \CER(G) \cap \bar Q_{w\rightarrow u}$. In the remainder of the proof, we show that
$\CER(G_{w\rightarrow u}) = \proj_{G_{w\rightarrow u}}(\overline{\CER}(G))$. This equivalence follows trivially for Cases 1 and 2, above. Henceforth, we prove it for Case~3.

Next, we characterize the convex hull of the set $\S^{=}(K_{\bar e})$ using RLT~\cite{SheAda90}. The product factors in this case are given by
$\psi(U, \bar e) = \prod_{v \in U} {(1-z_v)} \prod_{v \in \bar e \setminus U}{z_v}$ for all $U \subseteq \bar e$. In addition to imposing nonnegativity on all product factors, we should multiply the equality constraint $z_u = z_w$ by the product factors. Two cases arise:
\begin{itemize}

\item If $u \in U$ and $w \notin U$ or if $u \notin U$ and $w \in U$, multiplying $z_u = z_w$ by $\psi(U, \bar e)$, we obtain  $\psi(U,\bar e) = 0$.

\item If $u,w \notin U$ (resp. $u,w \in U$), then multiplying $z_u = z_w$ by $\psi(U,\bar e)$ yields $\psi(U,\bar e) = \psi(U,\bar e)$ (resp. $0 = 0$).
\end{itemize}
Therefore, the convex hull of $\S^{=}(K_{\bar e})$ is given by:
\begin{align}
& \psi(U, \bar e) = 0, \; \forall U \subseteq \bar e: |U\cap\{u,w\}| =1\label{RLTeq}\\ 
& \psi(U,\bar e) \geq 0, \; \forall U \subseteq \bar e: |U\cap\{u,w\}| = 0 \text{ or } 2. \label{RLTineq}
\end{align}
First, let $U=U_1=\{w\}$; substituting in~\eqref{RLTeq}, we get
$\psi(U_1,\bar e)= (1-z_w)\prod_{v\in \bar e \setminus \{w\}}{z_v} = 0$, which after linearization can be written as:
\begin{equation}\label{lm1}
    z_{\bar e} = z_{\bar e \setminus\{w\}}, \; \forall \bar e \in \ME: \bar e \supset \{u,w\}.
\end{equation}
Next, for any $v' \in \bar e \setminus \{u\}$, let $U=U_2=\{w,v'\}$;
substituting in~\eqref{RLTeq} we get $\psi(U_2, \bar e)= (1-z_w)(1-z_{v'})\prod_{v\in \bar e \setminus \{w,v'\}}{z_v} = (1-z_w)\prod_{v \in \bar e\setminus \{w,v'\}}{z_v}-\psi(U_1, \bar e)= 0$, which in turn implies that $\psi'(U_2,\bar e):=\psi(U_2,\bar e) + \psi(U_1,\bar e) =(1-z_w)\prod_{v \in \bar e\setminus \{w,v'\}}{z_v} = 0$, and this after linearization can be equivalently written as:
\begin{equation}\label{lm2}
    z_{\bar e \setminus \{v\}} = z_{\bar e \setminus\{w,v\}}, \; \forall v \in \bar e \setminus \{u\}, \; \forall \bar e \in \ME: \bar e \supset \{u,w\}.
\end{equation}
Next, for any $v' \neq v'' \in \bar e \setminus \{u\}$, let $U=U_3=\{w,v',v''\}$; substituting in~\eqref{RLTeq} we get $\psi(U_3, \bar e)= (1-z_w)(1-z_{v'})(1-z_{v''})\prod_{v\in \bar e \setminus \{w,v',v''\}}{z_v} = (1-z_w)\prod_{v\in \bar e \setminus \{w,v',v''\}}{z_v}-\psi'(U_2,\bar e)-\psi(U_2,\bar e) = (1-z_w)\prod_{v\in \bar e \setminus \{w,v',v''\}}{z_v}-\psi(U_1,\bar e)-2\psi(U_2, \bar e) = 0$. This in turn implies that $(1-z_w)\prod_{v\in \bar e \setminus \{w,v',v''\}}{z_v} = 0$, which after linearization can be equivalently written as:
\begin{equation}\label{lm3}
    z_{\bar e \setminus \{v,v'\}} = z_{\bar e \setminus\{w,v,v'\}}, \; \forall v \neq v'\in \bar e \setminus \{u\}, \; \forall \bar e \in \ME: \bar e \supset \{u,w\}.
\end{equation}
Therefore, by a recursive application of the above argument we find that for each $\bar e \in \ME$ with $\bar e \supset \{u,w\}$ we have:
\begin{equation}\label{lm4}
z_e = z_{e\setminus\{w\}}, \; \forall e \subseteq \bar e : e\supseteq \{u,w\}.
\end{equation}
Now let $U=U'_1=\{u\}$; substituting in~\eqref{RLTeq} and linearizing we get
$z_{\bar e} = z_{\bar e \setminus \{u\}}$, which together with~\eqref{lm1} implies that
$$
z_{\bar e \setminus \{u\}} = z_{\bar e \setminus \{w\}}, \; \forall \bar e \in \ME: \bar e \supset \{u,w\}.
$$
Similarly, replacing $w$ by $u$ in $U_2$ and $U_3$ and following a similar line of arguments and using~\eqref{lm4}, 
for each $\bar e \in \ME$ with $\bar e \supset \{u,w\}$ we get:
\begin{equation}\label{lm5}
    z_e = z_{(e\setminus\{w\})\cup \{u\}}, \; \forall e \subseteq \bar e: e \ni w, e \not\ni u. 
\end{equation}
Therefore, the convex hull of $\S^{=}(K_{\bar e})$ for any $\bar e \in \ME$ with $\bar e \supset \{u,w\}$ can be written as:
\begin{align}
    & z_{e} = z_{e\setminus \{w\}}, \quad \forall e \subseteq \bar e: \; |e| \geq 2, \; e\supset \{u, w\} \label{simpf3}\\
& z_{e} = z_{(e\setminus \{w\})\cup\{u\}}, \quad \forall e \subseteq \bar e: \; |e| \geq 2, \; e\ni w, e\not\ni u \label{simpf4}\\    
    & \psi(U, \bar e \setminus \{w\}) \geq 0, \quad \forall U \subseteq \bar e \setminus \{w\}\label{simpf5}.
\end{align}
By~\cref{prop:RLT}, inequalities~\eqref{simpf5} coincide with the description of $\MP(K_{\bar e \setminus \{w\}})$, where  $K_{\bar e \setminus \{w\}}$ is the complete hypergraph with node set $\bar e \setminus \{w\}$, which can be obtained from $K_{\bar e}$ by contracting $w$ to $u$.
By using equations~\eqref{simpf3} and~\eqref{simpf4} to project out variables $z_e$ for $e \subseteq \bar e$ with $e \ni w$, we conclude that $\CER(G_{w\rightarrow u}) = \proj_{G_{w\rightarrow u}}(\overline{\CER}(G))$. This together with $\overline{\CER}(G) \subseteq \CER(G)\cap \bar Q_{w\rightarrow u}$ implies that
$\CER(G_{w\rightarrow u}) \subseteq \proj_{G_{w\rightarrow u}}(\CER(G)\cap \bar Q_{w\rightarrow u})$.
\end{prf}



Thanks to~\cref{lemma:contract}, we can show that the tightness of the complete edge relaxation is preserved under the node contraction operation.

\begin{proposition}\label{cor: contraction}
Let $G=(V, E)$ be a hypergraph, let $u,w \in V$, and let $G_{w\rightarrow u}$ be the hypergraph obtained from $G$ by contracting $w$ to $u$.
    If $\CER(G)$ is an extension of $\MP(G)$, then $\CER(G_{w\rightarrow u})$ is an extension of $\MP(G_{w\rightarrow u})$.
\end{proposition}

\begin{prf}
    Denote by $\proj_G(\CER(G))$ the set obtained from $\CER(G)$
    by projecting out variables $z_e$, $e \in \cl(E) \setminus E$. By assumption, we have $\proj_G(\CER(G)) = \MP(G)$. Therefore,
    $\proj_G(\CER(G)) \cap Q_{w\rightarrow u} = \proj_G(\CER(G) \cap \bar Q_{w\rightarrow u}) = \MP(G) \cap Q_{w\rightarrow u}$, where the first equality follows since the operator $\proj_{G}(\cdot)$ does not project out variables $z_u, z_w$.
    We then deduce that
    \begin{align*}
    \proj_{G_{w\rightarrow u}}(\proj_G(\CER(G) \cap \bar Q_{w\rightarrow u})) = & \proj_{G}(\proj_{G_{w\rightarrow u}}(\CER(G) \cap \bar Q_{w\rightarrow u})) \\
    = &\proj_{G_{w\rightarrow u}}(\MP(G) \cap Q_{w\rightarrow u}).
    \end{align*}
    By Part~(i) of~\cref{lemma:contract} we have $\proj_{G_{w\rightarrow u}}(\MP(G) \cap Q_{w\rightarrow u}) = \MP(G_{w\rightarrow u})$, and by Part~(ii) of~\cref{lemma:contract} we have $\proj_{G_{w\rightarrow u}}(\CER(G) \cap \bar Q_{w\rightarrow u}) \supseteq \CER(G_{w\rightarrow u})$. These together imply that $\MP(G_{w\rightarrow u}) \supseteq \proj_G(\CER(G_{w\rightarrow u}))$. However, by definition $\MP(G_{w\rightarrow u}) \subseteq \proj_{G}(\CER(G_{w\rightarrow u}))$. Hence, we deduce that
$\MP(G_{w\rightarrow u}) = \proj_{G}(\CER(G_{w\rightarrow u}))$.
\end{prf}

\subsection{The node expansion operator}

Let $G=(V, E)$ be a hypergraph, let $w \in V$, and let $f$ be a nonempty set of nodes disjoint from $V$. 
We define the hypergraph obtained from $G$ by \emph{expanding $w$ to $f$} as the hypergraph $G'=\left(V', E'\right)$, where:
$$
\begin{aligned}
& V':=V \setminus \{w\} \cup f, \\
& E':=\{f\} \cup\{e \in E: w \notin e\} \cup\{e \setminus \{w\} \cup f : e \in E, w \in e\}.
\end{aligned}
$$
In the following we show how the multilinear polytope and the complete edge relaxation are transformed under the node expansion operator. To this end, we use the following result from \cite{dPDiG21IJO}, which allows us to combine perfect formulations that overlap in only one variable.

\begin{lemma}[lemma~1 in~\cite{dPDiG21IJO}]
\label{lem PQR}
Let $m_1,m_2,n_1,n_2$ be nonnegative integers with $m_1,m_2 \ge 1$.
Let $A \in \R^{m_1 \times (n_1 + 1)}$, 
$b \in \R^{m_1}$, 
$C \in \R^{m_2 \times (n_2 + 1)}$, 
$d \in \R^{m_2}$, 
and define
\begin{align*}
P := & \bra{(x, y) \in \R^{n_1} \times \R : A(x, y) \leq b}, \\
Q := & \bra{(y,z) \in \R \times \R^{n_2} : C(y,z) \leq d}.
\end{align*}
Assume that $P$ and $Q$ are polytopes with binary vertices.
Then the following polytope has binary vertices:
\begin{align*}
R = & 
\bra{(x, y, z) \in \R^{n_1} \times \R \times \R^{n_2} : A(x, y) \leq b, \ C(y,z) \leq d}.
\end{align*}
\end{lemma}

The next theorem shows how an extended formulation for the multilinear polytope is transformed under the node expansion operator and is of independent interest. The result is similar in nature to theorem 4 in \cite{dPDiG21IJO}. However, it concerns general extended formulations, instead of only formulations in the original space.
In the following, given a hypergraph $G=(V,E)$ and an edge $f \in E$, by the \emph{standard linearization of $f$}, we imply the system of inequalities in variables $y \in \R^{f \cup \bra{f}}$:
\begin{align}
\label{stdlin}
\begin{split}
y_f \geq \sum_{v\in f}{y_v}-|f|+1, \\
y_f \leq y_v & \qquad \forall v \in f, \\
y_f \geq 0, \\
y_v \leq 1 & \qquad \forall v \in f.
\end{split}
\end{align}

\begin{theorem}
\label{th not Silvia}
Let $G=(V, E)$ be a hypergraph, let $w \in V$, let $f$ be a nonempty set of nodes disjoint from $V$, and let $G'=(V', E')$ be the hypergraph obtained from $G$ by expanding $w$ to $f$.
Let $A(x,u) \le b$ be an extended formulation of $\MP(G)$, where $x\in \R^{V \cup E}$ are the original variables, and where $u \in \R^U$ denote the extended variables.
Let $C y \le d$ be the system \eqref{stdlin}, in variables $y \in \R^{f \cup \bra{f}}$.
Then, an extended formulation of the multilinear polytope $\MP(G')$, with original variables $z \in \R^{V' \cup E'}$ and extended variables $u \in \R^U$, is obtained from $A(x,u) \leq b$, $Cy \leq d$ by renaming the $x$ variables in $A(x,u) \leq b$ as follows:
\begin{equation}
\label{eq rename x}
\begin{aligned}
& x_w \to z_f, \\
& x_v \to z_v && \forall v \in V \setminus \bra{w}, \\
& x_e \to z_e && \forall e \in E, \ w \notin e, \\
& x_e \to z_{e \setminus \bra{w} \cup f} && \forall e \in E, \ w \in e,
\end{aligned}
\end{equation}
and by renaming the $y$ variables in $Cy \leq d$ as follows:
\begin{equation}
\label{eq rename y}
\begin{aligned}
& y_f \to z_f, \\
& y_v \to z_v && \forall v \in f.
\end{aligned}
\end{equation}
\end{theorem}

\begin{prf}
Let $V_x := V \setminus \bra{w}$, and observe that $V' = V_x \cup f$.
Let $E_x := \{e \in E: w \notin e\} \cup\{e \setminus \{w\} \cup f : e \in E, w \in e\}$, and observe that $E' = E_x \cup \bra{f}$.

Next, we consider the extended formulation of $\MP(G)$ given by the system $A(x,u) \le b$.
We rename the $x$ variables as in \eqref{eq rename x}, and obtain the system $A(z,u) \le b$ in variables $z \in \R^{V_x \cup E_x \cup \bra{f}}$ and $u \in \R^U$.
Let $P$ be the polyhedron defined by
\begin{align*}
P :=\bra{(z, u) \in \R^{V_x \cup E_x \cup \bra{f}} \times \R^{U} : A(z, u) \leq b}.
\end{align*}
Let $\tilde P$ be obtained from $P$ by projecting out variables $u \in \R^U$, i.e., 
\begin{align*}
\tilde P :=\bra{z \in \R^{V_x \cup E_x \cup \bra{f}} : \exists u \in \R^{U} \text{ such that } A(z, u) \leq b}.
\end{align*}
By assumption, $\tilde P$ is obtained from $\MP(G)$, in variables $x \in \R^{V \cup E}$, by renaming variables as in \eqref{eq rename x}.
In particular, $\tilde P$ is a polytope with binary vertices.
Denote by $\tilde Az \le \tilde b$ a system of linear inequalities, in variables $z \in \R^{V_x \cup E_x \cup \bra{f}}$, describing $\tilde P$, i.e.,
\begin{align*}
\tilde P =\bra{z \in \R^{V_x \cup E_x \cup \bra{f}} : \tilde Az \leq \tilde b}.
\end{align*}

Next, consider the standard linearization of the edge $f$, given by the system $Cy \le d$.
We rename the $y$ variables as in \eqref{eq rename y}, and obtain the system $Cz \le d$ in variables $z \in \R^{f \cup \bra{f}}$.
Let $Q$ be the polyhedron defined by
\begin{align*}
Q := \bra{z \in \R^{f \cup \bra{f}} : Cz \le d}.
\end{align*}
If we denote by $H$ the hypergraph $H = (f,\bra{f})$, 
it is well-know that $Q$ is obtained from $\MP(H)$, in variables $y \in \R^{f \cup \bra{f}}$, by renaming variables as in \eqref{eq rename y}.
In particular, $Q$ is a polytope with binary vertices.

Observe that $\tilde P$ and $Q$ share only one common variable, namely variable $z_f$.
In $\tilde P$, such variable was obtained by renaming variable $x_w$, while in $Q$ it was obtained by renaming variable $y_f$. 
We can then apply \cref{lem PQR} to $\tilde P$ and $Q$, and obtain that the following polytope has binary vertices:
\begin{align}
\label{eq use it later}
\begin{split}
R 
& = 
\bra{z \in \R^{V_x \cup f \cup E_x \cup \bra{f}} : \tilde Az \leq \tilde b, \ Cz \leq d} \\
& = 
\bra{z \in \R^{V' \cup E'} : \tilde Az \leq \tilde b, \ Cz \leq d}.
\end{split}
\end{align}

Recall that the system $\tilde Az \leq \tilde b$ is obtained from the system $A(z, u) \leq b$ by projecting out variables $u \in \R^U$.
Since the variables $u \in \R^U$ are not present in the system $Cz \leq d$, projecting out variables $u \in \R^U$ from the system $A(z, u) \leq b$, $Cz \leq d$, results in the system $\tilde Az \leq \tilde b$, $Cz \leq d$.
Therefore, $A(z, u) \leq b$, $Cz \leq d$ is an extended formulation of $R$.
More precisely,
\begin{align*}
R = & 
\bra{z \in \R^{V' \cup E'} : \exists u \in \R^U \text{ such that } A(z, u) \leq b, \ Cz \leq d}.
\end{align*}
To conclude the proof, we need to show that $R$ coincides with $\MP(G')$.
We already know that both $R$ and $\MP(G')$ have binary vertices, so we only need to check that they contain the same binary vectors.
The remainder of the proof is devoted to proving this claim.

Let $\bar{z}$ be a binary vector in $\MP(G')$.
We prove $\bar{z} \in R$.
From \eqref{eq use it later}, it suffices to show that the restriction of $\bar z$ to $\R^{V_x \cup E_x \cup \bra{f}}$ is in $\tilde P$, and the restriction of $\bar z$ to $\R^{f \cup \bra{f}}$ is in $Q$.
To do so, we define binary vectors $\bar x \in \R^{V \cup E}$ and $\bar y \in \R^{f \cup \bra{f}}$ using the inverse of \eqref{eq rename x} and \eqref{eq rename y}, and prove $\bar x \in \MP(G)$ and $\bar y \in \MP(H)$.
Specifically, $\bar x \in \R^{V \cup E}$ is defined by
\begin{equation}
\label{eq get x}
\begin{aligned}
& \bar{x}_w :=\bar{z}_f, \\
& \bar{x}_v :=\bar{z}_v && \forall v \in V \setminus \bra{w}, \\
& \bar{x}_e :=\bar{z}_e && \forall e \in E, \ w \notin e, \\
& \bar{x}_e :=\bar{z}_{e \setminus \bra{w} \cup f} && \forall e \in E, \ w \in e,
\end{aligned}
\end{equation}
and $\bar y \in \R^{f \cup \bra{f}}$ is defined by
\begin{equation}
\label{eq get y}
\begin{aligned}
& \bar{y}_f := \bar{z}_f, \\
& \bar{y}_v := \bar{z}_v && \forall v \in f.
\end{aligned}
\end{equation} 
We have $\bar y \in \MP(H)$ because
\begin{align*}
\bar{y}_f = \bar{z}_f = \prod_{v \in f} \bar{z}_v = \prod_{v \in f} \bar{y}_v.
\end{align*}
Furthermore, we have $\bar x \in \MP(G)$ because, for every $e \in E$ with $w \notin e$,
\begin{align*} 
\bar{x}_e & = \bar{z}_e 
= \prod_{v \in e} \bar{z}_v
= \prod_{v \in e} \bar{x}_v,
\end{align*}
and for every $e \in E$ with $w \in e$,
\begin{align*} 
\bar{x}_e & = \bar{z}_{e \setminus \bra{w} \cup f}
= \pare{\prod_{v \in e \setminus \bra{w}} \bar{z}_v} \pare{\prod_{v \in f} \bar{z}_v} \\
& = \pare{\prod_{v \in e \setminus \bra{w}} \bar{z}_v} \bar{z}_f 
= \pare{\prod_{v \in e \setminus \bra{w}} \bar{x}_v} \bar{x}_w
= \prod_{v \in e} \bar{x}_v.
\end{align*} 

Now, let $\bar{z}$ be a binary vector not in $\MP(G')$.
We prove $\bar{z} \notin R$.
From \eqref{eq use it later}, it suffices to show that either the restriction of $\bar z$ to $\R^{V_x \cup E_x \cup \bra{f}}$ is not in $\tilde P$, or the restriction of $\bar z$ to $\R^{f \cup \bra{f}}$ is not in $Q$.
To do so, we define binary vector $\bar x \in \R^{V \cup E}$ as in \eqref{eq get x}, and $\bar y \in \R^{f \cup \bra{f}}$ as in \eqref{eq get y}, and prove $\bar x \notin \MP(G)$ or $\bar y \notin \MP(H)$.
If $\bar{z}_f \neq \prod_{v \in f} \bar{z}_v$, then 
\begin{align*}
\bar{y}_f = \bar{z}_f \neq \prod_{v \in f} \bar{z}_v = \prod_{v \in f} \bar{y}_v,
\end{align*}
and $\bar y \notin \MP(H)$.
Thus, we now assume $\bar{z}_f = \prod_{v \in f} \bar{z}_v$.
Since $\bar{z}$ is not in $\MP(G')$, there exists an edge $g' \in E'$ with $g' \neq f$ such that $\bar{z}_{g'} \neq \prod_{v \in g'} \bar{z}_v$. 
Let $g$ be the edge of $G$ corresponding to $g'$, i.e., $g = g'$ if $g' \cap f = \emptyset$ and $g = g' \cup \bra{w} \setminus f$ if $f \subseteq g'$.
Note that in the first case we have $w \notin g$, while in the second case we have $w \in g$ and $g' = g \setminus \bra{w} \cup f$.
We now show $\bar{x} \notin \MP(G)$.
In fact, if $w \notin g$, then
\begin{align*} 
\bar{x}_g & = \bar{z}_{g'} 
\neq \prod_{v \in g'} \bar{z}_v
= \prod_{v \in g} \bar{x}_v.
\end{align*} 
On the other hand, if $w \in g$, then
\begin{align*} 
\bar{x}_g = \bar{z}_{g \setminus \bra{w} \cup f} = \bar{z}_{g'}
\neq \prod_{v \in g'} \bar{z}_v 
& =
\pare{\prod_{v \in g \setminus \bra{w}} \bar{z}_v} \pare{\prod_{v \in f} \bar{z}_v} = \\
& = \pare{\prod_{v \in g \setminus \bra{w}} \bar{z}_v} \bar{z}_f 
= \pare{\prod_{v \in g \setminus \bra{w}} \bar{x}_v} \bar{x}_w
= \prod_{v \in g} \bar{x}_v.
\end{align*} 
\end{prf}

Thanks to~\cref{th not Silvia}, we can show that the tightness of the complete edge relaxation is preserved under the node expansion operation.

\begin{proposition}\label{cor: expansion}
Let $G=(V, E)$ be a hypergraph, let $w \in V$, let $f$ be a nonempty set of nodes disjoint from $V$, and let $G'$ be the hypergraph obtained from $G$ by expanding $w$ to $f$.
If $\CER(G)$ is an extension of $\MP(G)$, then $\CER(G')$ is an extension of $\MP(G')$.
\end{proposition}

\begin{prf}
By assumption, $\CER(G)$ is an extension of $\MP(G)$. 
Denote by $x$ the original variables and by $u$ the extended variables in such extended formulation.
By~\cref{th not Silvia}, an extended formulation of $\MP(G')$, in original variables $z$ and extended variables $u$, is obtained from the description of $\CER(G)$ and \eqref{stdlin} by renaming the $x$ variables in the description of $\CER(G)$ as in \eqref{eq rename x} and renaming the $y$ variables in \eqref{stdlin} as in \eqref{eq rename y}. 
Let us consider the obtained formulation of $\MP(G')$. 
By construction, each inequality in this formulation contains variables corresponding to nodes and edges inside some maximal edge of $G'$. 
Moreover from the definition of complete edge relaxation, as well as \cref{prop:RLT}, the description of $\CER(G')$ is obtained by putting together the convex hull description for each subhypergraph of $G'$ induced by a maximal edge. 
It follows that all inequalities in the obtained extended formulation of $\MP(G')$ are implied by some inequalities in the description of $\CER(G')$. 
Therefore, $\CER(G')$ is an extension of $\MP(G')$. 
\end{prf}

\subsection{Reducing the length of simple cycles}

Thanks to the node contraction and node expansion operators introduced earlier, we are now able ``reduce'' the length of simple cycles. This in turn enables us to limit our attention to simple cycles of length three. The following theorem is of independent interest in hypergraph theory.

\begin{theorem}\label{lengthReduction}
Let $G = (V,E)$ be a hypergraph, assume that it contains a simple cycle, and let $\ell$ denote the minimum length of a simple cycle in $G$.
If $\ell \geq 4$, through a sequence of expansion and contraction operations, we can obtain, from $G$, a new hypergraph $G'$ that contains a simple cycle of length at most $\ell-1$.
\end{theorem}

\begin{prf}
Let $C=e_1,e_2,\ldots,e_{\ell},e_{\ell+1}$, for some $\ell \geq 4$, be a simple cycle of $G$ of minimum length.
As before, define $s_i := e_i \cap e_{i+1}$ for every $i \in [\ell]$. 
In the following, we define the expansion and contraction operations that enable us to reduce the length of a simple cycle.

\textbf{Expansion.}
Recursively, we expand each node $v \in s_1 \cup s_2$ to a set of new nodes $f_v$ of cardinality $\card{s_1 \cup s_2}$.
Since, for every node $v \in s_1 \cup s_2$, there is a bijection between the nodes in $f_v$ and the nodes in $s_1 \cup s_2$ (including $v$ itself), we can label each node in $f_v$ with the corresponding node in $s_1 \cup s_2$.

\textbf{Contraction.}
For each pair of distinct nodes $v,w \in s_1 \cup s_2$, we recursively contract the node in $f_v$ labeled by $w$ to the node in $f_w$ labeled by $v$;
we denote the resulting node by $u_{vw}$, where $vw$ is an unordered pair.
Note that each such contraction can be performed since each such pair of nodes, after expansion, is contained in $e_2$.
Denote by $G'=(V',E')$ the obtained hypergraph and by 
\begin{align*}
U & := \bra{u_{vw} : v,w \in s_1 \cup s_2, \ v \neq w}.
\end{align*}
For every node $v \in s_1 \cup s_2$, we identify the node in $f_v$ labeled by $v$ (which has not been contracted) with the original node $v$.
In this way, we have that $V'$ is the disjoint union of $V$ and $U$.

\textbf{The new simple cycle.}
For every $i = 1,\dots,\ell$, let $e'_i$ be the edge of $G'$ that originated from the edge $e_i$ of $G$.
More precisely, $e'_i$ is the edge of $G'$ given by

\begin{align}
\label{defnewej}
e'_i := e_i \cup \bra{u_{vw} \in U : e_i \cap \bra{v,w} \neq \emptyset}.
\end{align}

Define $C':=e'_1,e'_3,e'_4,\ldots,e'_{\ell},e'_{\ell+1}$.
Define $s'_1 := e'_1 \cap e'_3$ and, for every $i = 3,\dots,\ell$, define $s'_i := e'_i \cap e'_{i+1}$.
In the remainder of the proof, we show that $C'$ is a simple cycle of $G'$ of length $\ell-1$.
From \cref{def simple}, it suffices to show that $C'$ satisfies condition \eqref{eq simple}.

Let $i,j,l \in [\ell]\setminus \bra{2}$ with  $i < j < l$ and $e' \in E'$. 
We need to show $(s'_i \cup s'_j \cup s'_l) \setminus e' \neq \emptyset$.
For a contradiction, assume $(s'_i \cup s'_j \cup s'_l) \subseteq e'$.

First, assume that $e'$ is an edge of $G'$ contained in $U$.
In this case, we obtain $(s'_i \cup s'_j \cup s'_l) \subseteq e' \subseteq U$.
In particular, 
$$
s'_j \subseteq U.
$$
Note that $j \ge 3$. Applying \eqref{defnewej} to $e'_j$ and $e'_{j+1}$, we obtain
$$
s_j = s'_j \setminus U = \emptyset,
$$ 
which contradicts \cref{lem relate,obs simple}.

It remains to consider the case in which $e'$ is an edge of $G'$ that is not contained in $U$.
In this case, $e'$ is an edge that is originated from an edge of $G$.
More precisely, there exists $e \in E$ such that
\begin{align}
\label{defnewe}
e' := e 
\cup \bra{u_{vw} \in U : e \cap \bra{v,w} \neq \emptyset}.
\end{align}

\begin{claim}
\label{come on 1}
For every $k \in \bra{3,4,\dots,\ell}$, 
$s_k \setminus e \neq \emptyset$ implies $s'_k \setminus e' \neq \emptyset$.
\end{claim}

\begin{cpf}
By assumption, there exists $v \in s_k \setminus e$.
We show that $v \in s'_k \setminus e'$.

Applying \eqref{defnewej} to $e_k$ and $e_{k+1}$, we obtain $s_k \subseteq s'_k$, thus $v \in s_k$ implies $v \in s'_k$.
On the other hand, $v \in s_k$ implies $v \in V$ and $v \notin U$.
Since $v \notin e$, from \eqref{defnewe} we obtain $v \notin e'$.
\end{cpf}

\begin{claim}
\label{come on 2}
$s_1 \setminus e \neq \emptyset$ and $s_2 \setminus e \neq \emptyset$ imply $s'_1 \setminus e' \neq \emptyset$. 
\end{claim}

\begin{cpf}
By assumption, there exist $v \in s_1 \setminus e$ and $w \in s_2 \setminus e$.

Consider the case $v=w$.
In this case, we show $v \in s'_1 \setminus e'$.
First, observe that $v \in s'_1$.
In fact, by assumption, $v \in s_1 \cap s_2 = e_1 \cap e_2 \cap e_3$.
From \eqref{defnewej}, $v \in e'_1 \cap e'_2 \cap e'_3 \subseteq s'_1$.
On the other hand, from $v \in V$, $v \notin e$, and \eqref{defnewe}, we obtain $v \notin e'$.
This concludes the proof for the case $v=w$.

\smallskip

Next, consider the case $v \neq w$.
In this case, we show $u_{vw} \in s'_1 \setminus e'$.
First, observe that $u_{vw} \in s'_1$.
In fact, by assumption, $v \in s_1 = e_1 \cap e_2$, thus from \eqref{defnewej}, $u_{vw} \in e'_1 \cap e'_2$.
Similarly, by assumption, $w \in s_2 = e_2 \cap e_3$, thus from \eqref{defnewej}, $u_{vw} \in e'_2 \cap e'_3$.
Therefore, $u_{vw} \in e'_1 \cap e'_2 \cap e'_3 \subseteq s'_1$.
On the other hand, we have $u_{vw} \notin e'$.
In fact, from $v \notin e$, $w \notin e$, and \eqref{defnewe}, we obtain $u_{vw} \notin e'$.
This concludes the proof for the case $v \neq w$.
\end{cpf}

Recall that, in order to obtain a contradiction, we assumed $(s'_i \cup s'_j \cup s'_l) \subseteq e'$.
For $k \in \bra{i,j,l}$, we then have  $s'_k \subseteq e'$.
From \cref{come on 1}, if $k \in \bra{3,4,\dots, \ell}$, then $s_k \subseteq e$.
From \cref{come on 2}, if $k = 1$, then $s_1 \subseteq e$ or $s_2 \subseteq e$.
We have therefore obtained $p,q,r \in [\ell]$ with $p<q<r$ such that $(s_p \cup s_q \cup s_r) \subseteq e$, which gives us a contradiction.
\end{prf}

\subsection{The key structure}

By~\cref{lengthReduction} to prove~\cref{th necessity}, it essentially suffices to consider hypergraphs that contain simple cycles of length three. As we show in the next section, all such hypergraphs contain a subhypergraph with $n$ nodes for some $n\geq 3$, whose maximal edges consist of all subsets of cardinality $n-1$. In the following, we show that for such subhypergraphs, the complete edge relaxation is not an extension of the multilinear polytope. 

\begin{proposition}\label{almostfull}
    Let $G=(V,E)$ be a hypergraph; let $n :=|V| \geq 3$, and let the set of maximal edges of $G$ consist of all subsets of $V$ of cardinality $n-1$. Then $\CER(G)$ is not an extension of $\MP(G)$.
\end{proposition}

\begin{prf}
Denote by $\ME$ the set of maximal edges of $G$.
Consider the following inequality:
\begin{equation}\label{vmpg}
    \sum_{v \in V}{z_v} - \sum_{e \in \ME} {z_e} \leq |V|-2.
\end{equation}
We claim that~\eqref{vmpg} is a valid inequality for $\MP(G)$. To see this, notice that there are two cases in which this inequality can be violated:
\begin{itemize}
    \item $z_v = 1$ for all $v \in V$: in this case we also have $z_e = 1$ for all $e \in E$. Substituting this point in~\eqref{vmpg} we get $|V|-|\ME| = |V|-|V| \leq |V|-2$.

    \item $z_v = 1$ for all $v \in V' \subset V$, with $|V'|=|V|-1$: in this case we have $z_e = 1$ for some maximal edge $e \in \ME$, with $e = V'$. Substituting this point in~\eqref{vmpg} we get $|V|-1-1 \leq |V|-2$.
\end{itemize}
Hence, we conclude that~\eqref{vmpg} is a valid inequality for $\MP(G)$. Now consider the following fractional point $\tilde z \in \R^{V \cup \cl(E)}$:
\begin{equation}\label{fpoint}
    \tilde z_p = \frac{n-1-|p|}{n-1}, \quad \forall p \in V \cup \cl(E).
\end{equation}
First, consider the projection $\hat z$ of $\tilde z$ onto the space  $\R^{V \cup \ME}$; $\hat z_v = \frac{n-2}{n-1}$ for all $v \in V$ and $\hat z_e = 0$ for all $e \in \ME$. Substituting $\hat z$ in inequality~\eqref{vmpg} we obtain $n (\frac{n-2}{n-1})=n-\frac{n}{n-1} \not\leq n-2$, since $\frac{n}{n-1} \leq \frac{3}{2}$ for $n \geq 3$. Therefore, $\hat z$ violates the valid inequality~\eqref{vmpg}. In the following, we show that $\tilde z \in \CER(G)$ implying  $\CER(G)$ is not an extension of $\MP(G)$. The relaxation $\CER(G)$ is defined by the following inequalities:
\begin{equation}\label{rltstr}
\psi(U, e) \geq 0, \; \forall U\subseteq e, \; \forall e \in \ME,
\end{equation}
where $\psi(U, e)$ is defined by~\eqref{defRLT}. Define $m:=|U|$. Substituting~\eqref{fpoint} into $\psi(U, e)$, we obtain:
\begin{align*}
    \psi(U, e)  = &\frac{m}{n-1}-\binom{m}{1}\frac{m-1}{n-1}+\binom{m}{2}\frac{m-2}{n-1}-\ldots +(-1)^{m-1}\binom{m}{m-1}\frac{1}{n-1}+(-1)^m (0)\\
    =&\frac{1}{n-1}\Bigg(m-2\binom{m}{2}+3\binom{m}{3}-4\binom{m}{4}+\ldots+(-1)^{m-1}m \binom{m}{m}\Bigg)\\
    =& \frac{1}{n-1}\sum_{k=1}^m{(-1)^{k-1}k\binom{m}{k}}.
\end{align*}
where the second equality follows from the identity $\binom{n}{k}(n-k)=(k+1) \binom{n}{k+1}$. Moreover, we have
$$
\sum_{k=1}^m{(-1)^{k-1}k\binom{m}{k}}   
= \left\{ \begin{array}{ll}
        1 & \mbox{if $m = 1$};\\
        0 & \mbox{if $m = 0$ or $m > 1$}.
        \end{array} \right.
$$
Therefore, point~\eqref{fpoint} satisfies inequalities~\eqref{rltstr}, implying that $\tilde z \in \CER(G)$. We conclude that $\CER(G)$ is not an extension of $\MP(G)$. 
\end{prf}

Letting $n=3$ in~\cref{almostfull}, the hypergraph $G$ simplifies to a graph consisting of a cycle of length three.
Let $V=\{u,v,w\}$.
In this case inequality~\eqref{vmpg} simplifies to the following ``triangle'' inequality~\cite{Pad89}:
$$
z_u + z_v + z_w - z_{\{u,v\}}-z_{\{u,w\}}-z_{\{v,w\}} \leq 1.
$$
In fact, for any $n \geq 3$, the set of maximal edges $\ME$ of $G$ in~\cref{almostfull} forms a simple cycle. To see this, consider an ordering of the edges in $\ME$ denoted by $e_1, e_2, \ldots, e_n, e_{n+1}$, with $e_{n+1}=e_1$. We then have $s_i \cup s_j \cup s_k = V$, and by construction $V \setminus e \neq \emptyset$ for all $e \in E$.

\subsection{Proof of necessity}

We are now ready to prove~\cref{th necessity}.

\begin{prf}[Proof of \cref{th necessity}]

We prove that if $G$ is not $\alpha$-acyclic, then $\CER (G)$ is not an extension of $\MP(G)$. Suppose that $G$ is not $\alpha$-acyclic; then by~\cref{th characterization}, $G$ contains a simple cycle. Denote by $C= e_1, e_2, \ldots e_{\ell+1}$, where $e_{\ell+1} = e_1$, for some $\ell \geq 3$, a simple cycle of minimal length in $G$. Our proof is by induction on the length $\ell$ of $C$. Recall that $s_i = e_i \cap e_{i+1}$ for $i \in [\ell]$.

    \paragraph{Base case.} Suppose that $G$ contains a simple cycle of length three; $C= e_1, e_2, e_3, e_1$. 
    Define $\bar V := s_1 \cup s_2 \cup s_3$, and 
    let $L_{\bar V}$ denote the affine set as defined by~\eqref{affineL}. Denote by $G_{\bar V} = (\bar V, \bar E)$ the subhypergraph of $G$ induced by $\bar V$.
    By~\cref{cor: fixing}, to prove that $\CER(G)$ is not an extension of $\MP(G)$, it suffices prove that $\CER(G_{\bar V})$ is not an extension of $\MP(G_{\bar V})$. The hypergraph $G_{\bar V}$ is not $\alpha$-acyclic because it contains a simple cycle $C = \bar e_1, \bar e_2, \bar e_3, \bar e_{1}$, where $\bar e_i = e_i \cap \bar V$ for $i \in \{1,2,3\}$.
    In the following, for notational simplicity, instead of $\bar e_i$ we write $e_i$ for all $i \in [\ell+1]$.    
    We show that $\CER(G_{\bar V})$ is not an extension of $\MP(G_{\bar V})$. 
    Letting $\ell=3$ in condition~\eqref{eq simple}, we deduce that $\bar V$ is not an edge of $G_{\bar V}$, and that the following sets are nonempty:
    \begin{equation}\label{keynodes}
    \bar s_1 := (e_1 \cap e_2) \setminus e_3, \quad  
    \bar s_2 := (e_2 \cap e_3) \setminus e_1, \quad
    \bar s_3 :=(e_1 \cap e_3)\setminus e_2. 
    \end{equation}
    We consider two cases separately.

    \paragraph{Case 1.} 
    In this case, we assume that there exist three nodes $v_1 \in \bar s_1$, $v_2 \in \bar s_2$, $v_3 \in \bar s_3$ such that $\{v_1, v_2, v_3\}$ is not contained in an edge of $\bar E$.     
    Let $\tilde V = \{v_1, v_2, v_3\}$ and denote by $G_{\tilde V}$ the subhypergraph of $G_{\bar V}$ (or $G$) induced by $\tilde V$. By~\cref{cor: fixing}, it suffices to show that $\CER(G_{\tilde V})$ is not an extension of $\MP(G_{\tilde V})$.
    Clearly, $G_{\tilde V}$ is a graph consisting of a cycle of length three. By letting $n=3$ in~\cref{almostfull}, it follows that $\CER(G_{\tilde V})$ is not an extension of 
    $\MP(G_{\tilde V})$, which in turn implies that $\CER(G)$
    is not an extension of $\MP(G)$.

\paragraph{Case 2.} 
    In this case, we assume that for any three nodes $v_1 \in \bar s_1, v_2 \in \bar s_2, v_3 \in \bar s_3$, the set $\{v_1, v_2, v_3\}$ is contained in some edge of $\bar E$.     
    Recall that from condition~\eqref{eq simple} it follows that there is no edge $e \in \bar E$ such that $e = \bar V$. First, consider all subsets of $\bar V$ of cardinality $|\bar V|-1$. If all these subsets are edges of $G_{\bar V}$, then by~\cref{almostfull} we conclude that $\CER(G_{\bar V})$ is not an extension of $\MP(G_{\bar V})$. Otherwise, take a subset of $\bar V$ of cardinality $|\bar V|-1$, denoted by $\bar V^{-1}$,  which is not an edge of $G_{\bar V}$. Notice that $\bar V^{-1}$ contains at least one node from each of $\bar s_1, \bar s_2, \bar s_3$, since otherwise it coincides with one of the three edges $e_1$, $e_2$, $e_3$, implying that $\bar V^{-1}$ is an edge of $G_{\bar V}$. If all subsets of $\bar V^{-1}$ of cardinality $|\bar V|-2$ are contained in some edge of $\bar E$ and hence are edges of $G_{\bar V^{-1}}$, then by~\cref{almostfull} we conclude that $\CER(G_{\bar V^{-1}})$ is not an extension of $\MP(G_{\bar V^{-1}})$. Otherwise, we choose a subset of $\bar V^{-1}$ with cardinality $|\bar V|-2$ that is not contained in any edge of $\bar E$. We apply the above argument recursively, until we find a subset $\bar V^{-q}$ of cardinality $|\bar V| - q$ for some $q \in \{1,\ldots, |\bar V|-4\}$ satisfying the following properties:
    
\medskip
    \begin{itemize}
    \item [(i)] $\bar V^{-q}$ is not contained in any edge in $\bar E$, and 
    \item [(ii)]
    each subset of $\bar V^{-q}$ of cardinality $|\bar V^{-q}| -1$ is contained in some edge of $\bar E$. 
    \end{itemize}
\medskip

    We claim that we can always find a subset satisfying the two conditions above. To see this, consider any $V' \subseteq \bar V$ such that $|V'|=4$ and $V'\not\subseteq e$ for any  $e \in \bar E$. It then follows that $V'$ must contain three nodes $v_1 \in \bar s_1, v_2 \in \bar s_2, v_3 \in \bar s_3$. By assumption, we have $\{v_1, v_2, v_3\} \subseteq e$ for some $e \in \bar E$. 
    Denote by $\hat v$ the fourth node of $V'$.   Define $s_{\cap} = e_1 \cap e_2 \cap e_3$. We then have $\bar s_i = s_i \setminus s_{\cap}$, $i\in \{1,2,3\}$. Two cases arise:

    \medskip
    \begin{itemize}
    \item $\hat v \in \bar s_i$ for some $i \in \{1,2,3\}$: without loss of generality, let $\hat v \in \bar s_1$. Then by assumption, we have $\{\hat v, v_2, v_3\} \subseteq e$ for some $e \in \bar E$. Moreover, by construction,  
    $\{\hat v, v_1, v_2\} \subseteq e_2$ and
    $\{\hat v, v_1, v_3\} \subseteq e_1$.

    \item $\hat v \in s_{\cap}$: in this case, we have $\{\hat v, v_1, v_2\} \subseteq e_2$,
    $\{\hat v, v_1, v_3\} \subseteq e_1$, and
    $\{\hat v, v_2, v_3\} \subseteq e_3$.
    \end{itemize}
    \medskip

    Now consider the subhypergraph $G_{V'}$ of $G$ induced by $V'$. By assumption $V'$ is not an edge of $G_{V'}$.    
    By the above argument, all subsets of cardinality three of $V'$ are edges of $G_{V'}$. Hence, letting $n=4$ in~\cref{almostfull}, we deduce that
    $\MP(G_{V'}) \subset \proj_{G_{V'}}(\CER(G_{V'}))$. Therefore, we conclude that $\CER(G_{\bar V^{-q}})$ is not an extension of
    $\MP(G_{\bar V^{-q}})$, where $q \in \{1,\ldots,|\bar V|-4\}$. Hence, by~\cref{cor: fixing}, if $G$ contains a simple cycle of length three, $\CER(G)$ is not an extension of $\MP(G)$.
    
 \paragraph{The inductive step.} Now suppose that the minimal length of a simple cycle in $G$ is $\ell \geq 4$. Suppose that for any hypergraph $G$ that contains a simple cycle of length at most $\ell-1$, $\CER(G)$ is not an extension of $\MP(G)$. Then we would like to prove that for any hypergraph $G$ containing a simple cycle of length $\ell$,
 $\CER(G)$ is not an extension of $\MP(G)$ either.

 Let $C$ be a simple cycle of minimal length $\ell$ in $G$ for some $\ell \geq 4$. Then by~\cref{lengthReduction}, through a sequence of expansion and contraction operations, we can obtain a hypergraph $G'$ such that $G'$ contains a simple cycle of length at most $\ell-1$. By the induction hypothesis, $\CER(G')$ is not an extension of $\MP(G')$. Then by a repeated application of~\cref{cor: contraction} and~\cref{cor: expansion}, we deduce that $\CER(G)$ is not an extension of $\MP(G)$ and this completes the proof.
\end{prf}

\section{Generalized triangle inequalities}
\label{sec: inequalities}
In~\cite{Pad89}, Padberg introduced \emph{triangle inequalities}, a class of facet-defining inequalities for the Boolean quadric polytope. 
Let $G$ be a graph. Then for any cycle of length three in $G$ with nodes denoted by $u,v,w$, we define the following triangle inequalities:
\begin{subequations}
\label{trineq}
\begin{align}
       \label{trineq 1}
       & z_{\{u,v\}} + z_{\{u,w\}} \leq z_u + z_{\{v,w\}}\\
       \label{trineq 2}
       & z_{\{u,v\}} + z_{\{v,w\}} \leq z_v + z_{\{u,w\}} \\
       \label{trineq 3}
       & z_{\{u,w\}} + z_{\{v,w\}} \leq z_w + z_{\{u,v\}} \\
       \label{trineq 4}
       & z_u + z_v + z_w -z_{\{u,v\}}-z_{\{u,w\}} - z_{\{v,w\}} \leq 1.
\end{align}
\end{subequations}
Padberg proved that inequalities~\eqref{trineq} define facets of $\BQP(G)$.
He further generalized triangle inequalities to odd-cycle inequalities, a family of facet-defining inequalities corresponding to chordless cycles of any length. Subsequently, he showed that the polytope obtained by adding odd-cycle inequalities to the standard linearization coincides with the Boolean quadratic polytope if $G$ is a series-parallel graph.
In~\cite{BorCraHam92} the authors showed that triangle inequalities are Chv\'atal-Gomory (CG) cuts for the McCormick relaxation. Furthermore, if $G$ is a complete graph, they proved that the addition of triangle inequalities to the McCormick relaxation gives the CG-closure of $\BQPLP(G)$. 
On the computational front, triangle inequalities have had a significant impact on the performance of mixed-integer quadratic programming solvers~\cite{BonGunLin18}.

In this section, with the goal of strengthening the complete edge relaxation, we generalize triangle inequalities~\eqref{trineq} to obtain a family of facet-defining inequalities for simple cycles of length three. Notice that for cycles of length three, simple cycles coincide with $\alpha$-cycles. We show that the proposed inequalities are CG-cuts for the complete edge relaxation, while they are not CG cuts for the standard linearization. Throughout this section, we consider hypergraphs satisfying $G = \cl(G)$; \ie the space in which the complete edge relaxation is constructed.

\subsection{The switching operation for the multilinear polytope}
To define the new inequalities, we make use of the~\emph{switching operation} for the multilinear polytope, which has been employed in~\cite{BucRin08,dPKha17MOR}. 
This operation enables us to convert valid linear inequalities into other valid linear inequalities that
induce faces of the same dimension.
A similar operator has been introduced by several authors independently for Boolean quadric and cut polytopes ~\cite{Pad89,BarMah86}. 
Consider a hypergraph $G = (V,E)$ with $G =\cl(G)$.
Let $T \subseteq V$.
The characteristic vector for $T$, denoted by $\chi^T \in \bra{0,1}^{V \cup E}$, is defined by
\begin{equation*}
\begin{aligned}
& \chi^T_v :=
\begin{cases}
1 & \text{if $v \in T$} \\
0 & \text{if $v \notin T$}
\end{cases}
\qquad \forall v \in V, \\
& \chi^T_e :=
\begin{cases}
1 & \text{if $e \subseteq T$} \\
0 & \text{if $e \nsubseteq T$}
\end{cases}
\qquad \forall e \in E.
\end{aligned}
\end{equation*}
Clearly, the multilinear set of $G$ is precisely the set of all $\chi^T$, for $T \subseteq V$. 
Then, for any $U \subseteq V$, the corresponding \emph{switching} of $\chi^T$ is defined as 
\begin{align*}
\phi_U(\chi^T) := \chi^{T \Delta U}.
\end{align*}
The term `switching' naturally refers to flipping the value of each binary variable $z_v$, for every $v \in U$.
It is shown in theorem 4.3 in \cite{BucRin08} that, due to our assumption $G = \cl(G)$, each switching can be extended to an affine automorphism of $\MP(G)$ (and of $\R^{V \cup E}$) given by:
\begin{align}
\phi_U(z_v) & = \left \{
\begin{array}{lll}
1- z_v & \qquad & \qquad  \qquad  \text{if } v \in U\\
z_v & \qquad  & \qquad  \qquad \text{if } v \in V \setminus U
\end{array} \right. & \quad  \forall v \in V\label{mappx},\\
\phi_U(z_e) & =
\sum_{\substack{W \subseteq e \cap U \\ |W| \text{ even}}} {z_{(e \setminus U) \cup W}} - \sum_{\substack{W \subseteq e \cap U \\ |W| \text{ odd}}} {z_{(e \setminus U) \cup W}}
& \quad \forall e \in E, \label{mappy}
\end{align}
where we define $z_{\emptyset} := 1$.
%
For any $U \subseteq V$, we define the \emph{$U$-switching} of an inequality $a z \le \beta$ in $\R^{V \cup E}$, as the image of the inequality under $\phi_U$.
We also refer to a \emph{$U$-switching} of an inequality $a z \le \beta$, for some $U \subseteq V$, simply as a \emph{switching} of $a z \le \beta$.
Notice that if we set $U = \emptyset$, the $U$-switching of $a z \leq \beta$ coincides with $a z \leq \beta$.
It then follows from the above discussion (see corollary 4.5 in \cite{BucRin08}) that an inequality $a z \le \beta$ is valid for $\MP(G)$ if and only if any of its switchings is valid for $\MP(G)$, and the former is facet-defining for $\MP(G)$ if and only if the latter is.
We also refer to \emph{all switchings} of $\alpha z \leq \beta$ as the system of $U$-switchings of $a z \le \beta$, for all $U \subseteq V$.

Observe that the transformation defined by ~\eqref{mappx} and~\eqref{mappy} is also an affine automorphism of the complete edge relaxation $\CER(G)$.
In fact, the system of inequalities defining $\CER(G)$ is precisely the system obtained by putting together all switchings of the inequality $z_e \ge 0$, for every $e \in E_{\max}$.
This implies that an inequality $a z \le \beta$ is valid for $\CER(G)$ if and only if any of its switchings is valid for $\CER(G)$.
Since switchings preserve the integrality of the coefficients and the integrality of the right-hand side, we deduce that an inequality $a z \le \beta$ is a CG cut for $\CER(G)$ if and only if any of its switchings is a CG-cut for $\CER(G)$.


\begin{observation}\label{switchTri}
Let us consider Padberg's triangle inequalities for the Boolean quadric polytope.
It can be checked that every inequality in~\eqref{trineq} is a switching of inequality \eqref{trineq 1}:
Inequality~\eqref{trineq 2} can be obtained with $U=\{w\}$, 
inequality~\eqref{trineq 3} with $U=\{v\}$,
and inequality~\eqref{trineq 4} with $U=\{u\}$.
In fact, system~\eqref{trineq} consists of all switchings of any of its inequalities.
\end{observation}

\subsection{The new inequalities}

In the following, given an $\alpha$-cycle $C=e_1, e_2, e_3, e_1$, we define the support hypergraph of $C$ as the hypergraph $G_C=(V_C, E_C)$, where 
\begin{equation}\label{nec}
V_C := \cup_{i\in  [\ell]}{e_i}, \quad E_C:=\big\{e \subseteq e_i: |e| \geq 2, \forall i \in [\ell]\big\}.
\end{equation} 
We are now ready to present our new inequalities, which serve as the generalization of triangle inequalities, to the multilinear polytope.

\begin{proposition}\label{prop: gtri}
    Let $G$ be a hypergraph with $G= \cl(G)$, and let $C=e_1, e_2, e_3, e_1$ be an $\alpha$-cycle of $G$.  Suppose that
        \begin{equation}\label{wlog}
        e_1 \setminus (e_2 \cup e_3) = \emptyset, \;\;
        e_2 \setminus (e_1 \cup e_3) = \emptyset, \;\;
        e_3 \setminus (e_1 \cup e_2) = \emptyset.
    \end{equation}
    Then all switchings of the following inequalities are valid for $\MP(G)$:
\begin{subequations}
\label{gtrineq}
    \begin{align}
       \label{gtrineq 1}
       & z_{e_1} + z_{e_2} \leq z_{e_1 \cap e_2} + z_{e_3}\\
       \label{gtrineq 2}
       & z_{e_1} + z_{e_3} \leq z_{e_1 \cap e_3} + z_{e_2} \\
       \label{gtrineq 3}
       & z_{e_2} + z_{e_3} \leq z_{e_2 \cap e_3} + z_{e_1} \\
       \label{gtrineq 4}
       & z_{e_1 \cap e_2} + z_{e_1 \cap e_3} + z_{e_2 \cap e_3} -z_{e_1}-z_{e_2} - z_{e_3} \leq z_{e_1 \cap e_2 \cap e_3},
\end{align}
\end{subequations}
where we define $z_{\emptyset} := 1$. 
Moreover, all switchings of inequalities~\eqref{gtrineq} define facets of $\MP(G_C)$.
\end{proposition}

\begin{prf}
First, we show that it suffices to consider inequality \eqref{gtrineq 1} by proving that every other inequality in the system \eqref{gtrineq} is a switching of \eqref{gtrineq 1}. 
Let $U=e_1 \setminus e_2$; then the $U$-switching of \eqref{gtrineq 1} is given by 
$$
z_{e_1 \cap e_2}-z_{e_1} + z_{e_2} \leq z_{e_1 \cap e_2} + z_{e_3 \setminus (e_1 \setminus e_2)}-z_{e_3}.
$$
First note that 
$e_3 \setminus (e_1 \setminus e_2) = (e_2 \cap e_3) \cup (e_3 \setminus e_1)$. Moreover, by assumption~\eqref{wlog} we have $e_3 \setminus (e_1 \cup e_2) = \emptyset$, implying that $(e_2 \cap e_3) \cup (e_3 \setminus e_1)= e_2 \cap e_3$.
Substituting into the above inequality we get $z_{e_2} + z_{e_3} \leq z_{e_2 \cap e_3}+z_{e_1}$, \ie inequality~\eqref{gtrineq 3}. 
Similarly, we can prove that inequality~\eqref{gtrineq 2} is the $U$-switching of inequality~\eqref{gtrineq 1} with $U= e_2 \setminus e_1$. 
Next, let $U= (e_1 \cap e_2) \setminus e_3$;
then the $U$-switching of inequality~\eqref{gtrineq 1} is given by:
$$
z_{e_1 \setminus (e_2 \setminus e_3)} - z_{e_1}+z_{e_2 \setminus (e_1 \setminus e_3)} - z_{e_2} \leq z_{e_1 \cap e_2 \cap e_3}-z_{e_1 \cap e_2}+ z_{e_3}.
$$
Again, from assumption~\eqref{wlog} it follows that $e_1 \setminus (e_2 \setminus e_3) = e_1 \cap e_3$ and $e_2 \setminus (e_1 \setminus e_3) = e_2 \cap e_3$. 
Substituting into the above inequality yields inequality~\eqref{gtrineq 4}.

Henceforth, in the remainder of the proof, we consider inequality~\eqref{gtrineq 1}:
\begin{equation}\label{ftr}
    -z_{e_1 \cap e_2}+ z_{e_1} + z_{e_2} - z_{e_3}\leq 0.  
\end{equation}
We start by proving the validity of inequality~\eqref{ftr} for $\MP(G)$. Two cases arise:
\medskip
\begin{itemize}
    \item if $z_{e_1} = 1$ and $z_{e_2} = 0$ (or $z_{e_1} = 0$ and $z_{e_2} = 1$), then we have $z_{e_1 \cap e_2} = 1$, implying the validity of~\eqref{ftr}.
    \item if $z_{e_1} = z_{e_2} = 1$, then we have $z_{e_1 \cap e_2} = 1$, and by assumption~\eqref{wlog} we have $z_{e_3}= 1$, implying the validity of~\eqref{ftr}.
\end{itemize}
\medskip
Therefore, inequality~\eqref{ftr} is valid for $\MP(G)$. 

\medskip
Next, we show that inequality~\eqref{ftr} defines a facet of $\MP(G_C)$, where $G_C=(V_C, E_C)$ and $V_C, E_C$ are defined by~\eqref{nec}.
To this end, denote by $a z \le \alpha $ a nontrivial valid inequality for $\MP(G_C)$ that is satisfied tightly by all binary points in $\MP(G_C)$ that are binding for inequality~\eqref{ftr}. We show that the two inequalities coincide up to a positive scaling, which by the full dimensionality of $\MP(G_C)$~\cite{dPKha17MOR} implies that inequality~\eqref{ftr} defines a facet of $\MP(G_C)$.
Define $\bar s_1 := (e_1 \cap e_2) \setminus e_3$,  $\bar s_2 := (e_2 \cap e_3) \setminus e_1$, $\bar s_3 :=(e_1 \cap e_3)\setminus e_2$, and $s_\cap = e_1 \cap e_2 \cap e_3$. From assumption~\eqref{wlog} it follows that $e_i = \bar s_i \cap s_{\cap}$, $i \in \{1,2,3\}$. Recall that from condition~\eqref{eq simple} it follows that $\bar s_i \neq \emptyset$ for $i \in \{1,2,3\}$. In the following, by a \emph{tight point}, we imply a binary point in $\MP(G_C)$ that satisfies inequality~\eqref{ftr} tightly.

First, consider a tight point with $z_v = 0$ for all $v \in V_C$. Substituting into $a z \leq \alpha$ yields:
\begin{equation}\label{first}
   \alpha = 0.    
\end{equation}
Next, consider a tight point with $z_{\tilde v} = 1$ for some $\tilde v \in \bar s_2 \cup \bar s_3 \cup s_{\cap}$ and $z_v = 0$ otherwise. Substituting into $a z \leq \alpha$ and using~\eqref{first} and the fact that $\bar s_1 \neq \emptyset$, we obtain:
\begin{equation}\label{second}
   a_v = 0, \quad \forall v \in \bar s_2 \cup \bar s_3 \cup s_{\cap}. 
\end{equation}
If $|e_1 \cap e_2| > 1$, consider a tight point with $z_{\tilde v} = 1$ for some $\tilde v \in \bar s_1$ and $z_v = 0$ otherwise. Substituting into $a z \leq \alpha$ and using~\eqref{first}, we obtain:
\begin{equation}\label{third}
   a_v = 0, \quad \forall v \in \bar s_1, \text{ if } |e_1 \cap e_2| > 1. 
\end{equation}
If $|e_3| > 2$, let $\tilde e \subset e_3$ be an edge of cardinality two. Consider a tight point with $z_v =1$ for all $v \in \tilde e$ and $z_v = 0$ otherwise. Substituting into $a z \leq \alpha$ and using~\eqref{first} and~\eqref{second}, we obtain $a_{\tilde e} = 0$. If $|e_3| > 3$, let $\tilde e \subset e_3$ be an edge of cardinality three.  Again, consider a tight point as constructed above. Substituting into $a z \leq \alpha$ and using the above result, we get $a_{\tilde e} = 0$. By a recursive application of this argument, we get
\begin{equation}\label{fourth}
    a_e = 0, \quad \forall e \in E_C: e \subsetneq e_3. 
\end{equation}
If $|e_1| > 2$ (resp. $|e_2| > 2$), let $\tilde e \subset e_1$ (resp. $\tilde e \subset e_2$) be an edge of cardinality two such that $\tilde e \not\supseteq e_1 \cap e_2$. Consider a tight point with $z_v =1$ for all $v \in \tilde e$ and $z_v = 0$ otherwise. Substituting into $a z \leq \alpha$ and using~\eqref{first} and~\eqref{third} we obtain $a_{\tilde e} = 0$. Next we let $\tilde e$ to an edge of cardinality three such that $\tilde e \not\supseteq e_1 \cap e_2$. Applying this argument recursively, we obtain:
\begin{equation}\label{fifth}
    a_e = 0, \quad \forall e \in E_C: e \subsetneq e_1 \text{ or } e \subsetneq e_2: e \not\supseteq e_1 \cap e_2. 
\end{equation}
If $|\bar s_3| > 1$, consider a tight point with $z_v = 1$ for all $v \in e_2$ and $z_v = 0$, otherwise. Substituting into $a z \leq \alpha$ and using~\eqref{first}--\eqref{third} and~\eqref{fifth} we get:
\begin{equation}\label{sixth}
    \sum_{\substack{e \in E_C: \; e \subseteq e_2, \\
    e \supseteq e_1 \cap e_2}}{a_e} = 0.
\end{equation}
Now, let $\tilde v \in \bar s_3$.  Consider a tight point with  $z_v = 1$ for all $v \in e_2 \cup \{\tilde v\}$ and $z_v = 0$, otherwise. Substituting into $a z \leq \alpha$ and using~\eqref{first}--\eqref{fifth} we obtain:
$$
    \sum_{\substack{e \in E_C: e \subseteq e_2, \\
    e \supseteq e_1 \cap e_2}}{a_e} + a_{(e_1 \cap e_2) \cup\{\tilde v\}}= 0,
$$
which together with~\eqref{sixth} implies that  $a_{(e_1 \cap e_2) \cup\{\tilde v\}} = 0$. Next we should let $\tilde v, \bar v \in \bar s_3$ and construct a tight point as above to get $a_{(e_1 \cap e_2) \cup\{\tilde v, \bar v\}} = 0$. Using this argument recursively together with a symmetric argument, we deduce that:
\begin{equation}\label{seven}
    a_e = 0, \quad \forall  e \subsetneq e_1 \text{ or } e \subsetneq e_2: e \supset e_1 \cap e_2.
\end{equation}
Next, consider a tight point $z_v = 1$ for all $v \in e_1$ and $z_v = 0$, otherwise.  Substituting into $a z \leq \alpha$ and using~\eqref{first}--\eqref{seven} we obtain:
\begin{equation}\label{eight}
    a_{e_1 \cap e_2}  + a_{e_1} = 0.
\end{equation}
Moreover, by symmetry we have:
\begin{equation}\label{nine}
    a_{e_1 \cap e_2} + a_{e_2} = 0.
\end{equation}
Consider a tight point with $z_v = 1$ for all $v \in V_C$. Substituting into $a z \leq \alpha$ and using~\eqref{first}--\eqref{seven} we obtain:
\begin{equation}\label{ten}
    a_{e_1 \cap e_2}  + a_{e_1} + a_{e_2} + a_{e_3} = 0.
\end{equation}
Finally, consider a point $z_v = 1$ for all $v \in e_3$ and $z_v = 0$ otherwise. Since by assumption $\bar s_1 \neq \emptyset$, substituting this point in $az \leq \alpha$ and using~\eqref{first}--\eqref{ten} we get $a_{e_3} \leq 0$.
Since $az \leq \alpha$ is a nontrivial valid inequality for 
$\MP(G_C)$, from~\eqref{first}--\eqref{ten}, we deduce that this inequality can be equivalently written as $\beta (-z_{e_1 \cap e_2}+z_{e_1}+z_{e_2}-z_{e_3}) \leq 0$ for some $\beta > 0$. Therefore, inequality~\eqref{ftr} defines a facet of $\MP(G_C)$ and this completes the proof.
\end{prf}

\begin{observation}
    Let us consider assumption~\eqref{wlog} of~\cref{prop: gtri}. It is simple to check that if this assumption is not satisfied for an $\alpha$-cycle $C$ in $G$, then inequalities~\eqref{gtrineq} are not valid for $\MP(G)$. Now, let $G=(V,E)$ be a hypergraph with $G=\cl(G)$, and let  $C=e_1, e_2, e_3, e_1$ be an $\alpha$-cycle of $G$ that does not satisfy assumption~\eqref{wlog}. Define $\bar e_1 = e_1 \cap (e_2 \cup e_3)$, $\bar e_2 = e_2 \cap (e_1 \cup e_3)$, and $\bar e_3 = e_3 \cap (e_1 \cup e_2)$. Notice that since by assumption $G= \cl(G)$, we have $\bar e_1, \bar e_2, \bar e_3 \in E$. Moreover, since $e_1 \cap e_2 = \bar e_1 \cap \bar e_2$, $e_1 \cap e_3 = \bar e_1 \cap \bar e_3$, and $e_2 \cap e_3 = \bar e_2 \cap \bar e_3$, 
    from condition~\eqref{eq simple} it follows that $\bar C = \bar e_1, \bar e_2, \bar e_3, \bar e_1$ is an $\alpha$-cycle of $G$. Moreover,  $\bar C$ satisfies assumption~\eqref{wlog} and therefore we can write inequalities~\eqref{gtrineq} for $\bar C$ instead.   
\end{observation}


As we described in~\cref{switchTri}, in the case of triangle inequalities, the system~\eqref{trineq} represents all switchings of any of its inequalities. However, this is not the case for the generalized triangle inequalities~\eqref{gtrineq}.
The following example demonstrates how to obtain new switchings of inequalities in~\eqref{gtrineq}.

\begin{example}
Let $G$ be a hypergraph and suppose that $C=e_1, e_2, e_3, e_1$
with $e_1 = \{1,2,4\}$, $e_2=\{2,3,4\}$, $e_3=\{1,3,4\}$ is an $\alpha$-cycle of $G$.  Inequalities~\eqref{gtrineq} in this case are given by:
\begin{align*}
    \begin{split}
       & z_{\{1,2,4\}} + z_{\{2,3,4\}} \leq z_{\{2,4\}} + z_{\{1,3,4\}}\\
       & z_{\{1,2,4\}} + z_{\{1,3,4\}} \leq z_{\{1,4\}} + z_{\{2,3,4\}} \\
       & z_{\{2,3,4\}} + z_{\{1,3,4\}} \leq z_{\{3,4\}} + z_{\{1,2,4\}} \\
       & z_{\{2,4\}} + z_{\{1,4\}} + z_{\{3,4\}} -z_{\{1,2,4\}}-z_{\{2,3,4\}} - z_{\{1,3,4\}} \leq z_{4}.
    \end{split}
\end{align*}
Now, let $U=\{4\}$. Then the image of the above inequalities under the mapping $\phi_U$ defined by~\eqref{mappx} and~\eqref{mappy}, is given by:
\begin{align*}
    \begin{split}
       & z_{\{1,2\}}-z_{\{1,2,4\}} + z_{\{2,3\}}-z_{\{2,3,4\}} \leq z_2 - z_{\{2,4\}} + z_{\{1,3\}}- z_{\{1,3,4\}}\\
       & z_{\{1,2\}}-z_{\{1,2,4\}}+z_{\{1,3\}}- z_{\{1,3,4\}} \leq z_1 - z_{\{1,4\}} + z_{\{2,3\}}-z_{\{2,3,4\}} \\
       & z_{\{2,3\}}- z_{\{2,3,4\}} + z_{\{1,3\}}-z_{\{1,3,4\}} \leq z_3-z_{\{3,4\}} +z_{\{1,2\}}- z_{\{1,2,4\}} \\
       & z_1+z_2+z_3+z_4-z_{\{1,2\}}- z_{\{1,3\}}-z_{\{1,4\}} -z_{\{2,3\}}-z_{\{2,4\}} -
       z_{\{3,4\}}+z_{\{1,2,4\}}+z_{\{1,3,4\}} +z_{\{2,3,4\}}  \leq 1.
    \end{split}
\end{align*}
Next, suppose that $C=e_1, e_2, e_3, e_1$
with $e_1 = \{1,2,3\}$, $e_2=\{1,2,4\}$, $e_3=\{3,4\}$ is an $\alpha$-cycle of $G$. Inequalities~\eqref{gtrineq} in this case are given by:
\begin{align*}
    \begin{split}
       & z_{\{1,2,3\}} + z_{\{1,2,4\}} \leq z_{\{1,2\}} + z_{\{3,4\}}\\
       & z_{\{1,2,3\}} + z_{\{3,4\}} \leq z_{3} + z_{\{1,2,4\}} \\
       & z_{\{1,2,4\}} + z_{\{3,4\}} \leq z_{4} + z_{\{1,2,3\}} \\
       & z_{\{1,2\}} + z_{3} + z_{4} -z_{\{1,2,3\}}-z_{\{1,2,4\}} - z_{\{3,4\}} \leq 1.
    \end{split}
\end{align*}
Now, let $U=\{2\}$. Then the image of the above inequalities under the mapping $\phi_U$ defined by~\eqref{mappx} and~\eqref{mappy}, is given by:
\begin{align*}
    \begin{split}
       & z_{\{1,3\}}-z_{\{1,2,3\}} + z_{\{1,4\}}-z_{\{1,2,4\}} \leq z_1-z_{\{1,2\}} + z_{\{3,4\}}\\
       & z_{\{1,3\}}-z_{\{1,2,3\}} + z_{\{3,4\}} \leq z_{3} + z_{\{1,4\}}-z_{\{1,2,4\}} \\
       & z_{\{1,4\}}-z_{\{1,2,4\}} + z_{\{3,4\}} \leq z_{4} + z_{\{1,3\}}-z_{\{1,2,3\}} \\
       & z_1+ z_{3} + z_{4} -z_{\{1,2\}} -z_{\{1,3\}}-z_{\{1,4\}}- z_{\{3,4\}}+z_{\{1,2,3\}}+z_{\{1,2,4\}}  \leq 1.
    \end{split}
\end{align*}
$\diamond$
\end{example}

We leave open the question of whether the generalized triangle inequalities, together with the complete edge relaxation, characterize the multilinear polytope of the support hypergraph of an $\alpha$-cycle of length three.

\subsection{Generalized triangle inequalities as CG-cuts}

In the following, we show that inequalities~\eqref{gtrineq} can be obtained from the complete edge relaxation as CG-cuts. 
Together with~\cref{{prop: gtri}}, this demonstrates the usefulness of the complete edge relaxation, as it can be used to obtain new families of facet-defining inequalities for the multilinear polytope of general hypergraphs. Interestingly, inequalities~\eqref{gtrineq} are not CG-cuts for the standard linearization, another reason that conveys the superiority of the complete edge relaxation.

\begin{proposition}\label{cgcuts}
    Let $G=(V,E)$ be a hypergraph with $G= \cl(G)$, and let $C=e_1, e_2, e_3, e_1$ be an $\alpha$-cycle of $G$ satisfying assumption~\eqref{wlog}. Then all switchings of inequalities~\eqref{gtrineq} are CG-cuts for the complete edge relaxation $\CER(G)$.
\end{proposition}

\begin{prf}
From our discussion at the beginning of \cref{sec: inequalities}, 
it suffices to prove the statement for one of the inequalities in~\eqref{gtrineq}, for example for inequality~\eqref{gtrineq 1}.
%

Recall that we denote by $G_C=(V_C, E_C)$ the support hypergraph of the $\alpha$-cycle $C$, where $V_C$ and $E_C$ are defined by~\eqref{nec}. It is simple to check that inequalities defining $\CER(G_C)$ are implied by the inequalities defining $\CER(G)$.
To complete the proof, it suffices to show that the inequality
\begin{equation}\label{goal}
    2 z_{e_1 \cap e_2} - 2 z_{e_1}- 2 z_{e_2}+2 z_{e_3} \geq -1
\end{equation}
is implied by the inequalities defining $\CER(G_C)$, since inequality~\eqref{gtrineq 1} can then be obtained by dividing inequality~\eqref{goal} by 2 and rounding up the right-hand side.
To this end, we make use of the following two claims:
\begin{claim}\label{part1}
    The following inequalites are obtained by summing up some of the inequalities defining $\CER(G_C)$:
    \begin{align}
    \begin{split}\label{firstset}
        & z_{e_1 \cap e_2} - z_{e_1} \geq 0, \quad z_{e_1 \cap e_2} - z_{e_2} \geq 0\\
        & z_{e_1 \cap e_3} - z_{e_1} \geq 0, \quad z_{e_2 \cap e_3} - z_{e_2} \geq 0.
        \end{split}
    \end{align}
\end{claim}
\begin{cpf}
     By symmetry it suffices to show the statement for $z_{e_1 \cap e_2} - z_{e_1} \geq 0$. Consider the following inequalities defining $\CER(G_C)$:
     $$
     \psi\big(e_1 \setminus (e_2\cup J), e_1\big) \geq 0, \quad \forall J \subsetneq e_1 \setminus e_2.
     $$
     Summing up all of the above inequalities and using $z_e = \prod_{v \in e}{z_v}$ for all $e \in E$ we obtain:
     \begin{align*}
         \sum_{J \subsetneq e_1 \setminus e_2}{\psi\big(e_1 \setminus (e_2\cup J), e_1\big)}=&\prod_{v\in e_1 \cap e_2}{z_v}\sum_{J \subsetneq e_1 \setminus e_2}\Big(\prod_{v \in J}{z_v}\prod_{v \in e_1 \setminus (e_2 \cup J)}{(1-z_v)}\Big)\\
         =& \prod_{v\in e_1 \cap e_2}{z_v}\Big(1-\prod_{v \in e_1 \setminus e_2}{z_v}\Big) \\
         =& z_{e_1 \cap e_2}-z_{e_1} \geq 0.
     \end{align*}
\end{cpf}

\begin{claim}\label{part2}
    The following inequality is obtained by summing up some of the inequalities defining $\CER(G_C)$:
    \begin{equation}\label{secondset}
        1-z_{e_1 \cap e_3}-z_{e_2 \cap e_3}+z_{e_3} \geq 0.
    \end{equation}
\end{claim}
\begin{cpf}
Consider the following inequalities defining $\CER(G_C)$:
     $$
     \psi(e_3 \setminus J, e_3) \geq 0, \quad \forall J \subset e_3: J \not\supseteq e_1 \cap e_3, \; J \not\supseteq e_2 \cap e_3.
     $$
     Summing up all of the above inequalities and using $z_e = \prod_{v \in e}{z_v}$ for all $e \in E$ we obtain:
     \begin{align*}
         \sum_{\substack{J \subseteq e_3: \\ J \not\supseteq e_1 \cap e_3, \; J \not\supseteq e_2 \cap e_3}}{\psi(e_3 \setminus J, e_3)}=&1-\sum_{\substack{J \subseteq e_3: \\ J \supseteq e_1 \cap e_3}}{\psi(e_3 \setminus J, e_3)}-\sum_{\substack{J \subseteq e_3: \\ J \supseteq e_2 \cap e_3}}{\psi(e_3 \setminus J, e_3)}+\psi(\emptyset, e_3)\\
         =&1-\prod_{v \in e_1 \cap e_3}{z_v} \sum_{J \subseteq e_3 \setminus e_1}\Big(\prod_{v\in J}{z_v}\prod_{v\in e_3 \setminus (e_1 \cup J)}{(1-z_v)}\Big)\\
         &-\prod_{v \in e_2 \cap e_3}{z_v} \sum_{J \subseteq e_3 \setminus e_2}\Big(\prod_{v\in J}{z_v}\prod_{v\in e_3 \setminus (e_2 \cup J)}{(1-z_v)}\Big)+\prod_{v\in e_3}{z_v}\\
         =&1-z_{e_1 \cap e_3}-z_{e_2 \cap e_3} + z_{e_3} \geq 0,
     \end{align*}
     where the first equality follows from assumption~\eqref{wlog} and the second equality follows from the identity $\sum_{J \subseteq N}{\prod_{v\in J}{z_v}\prod_{N \setminus J}(1-z_v)}=1$ for any $N \subseteq V$.
\end{cpf}
Finally, summing up all of inequalities~\eqref{firstset} and inequality~\eqref{secondset} and inequality $z_{e_3} \geq 0$, we obtain inequality~\eqref{goal}.
\end{prf}

A natural question is whether inequalities~\eqref{gtrineq} are also CG-cuts for the standard linearization. The following example shows that the answer is negative.

\begin{example}
Let $G=(V,E)$ be a hypergraph with $G= \cl(G)$ and the set of maximal edges $E_{\max}= \{e_1, e_2, e_3\}$, where $e_1 = \{1,2,4\}$, $e_2=\{2,3,4\}$, $e_3=\{1,3,4\}$. It can be checked that $C=e_1, e_2, e_3, e_1$ is an $\alpha$-cycle of $G$ satisfying~\eqref{wlog}.  A generalized triangle inequality in this case is given by:
$$
z_{\{1,2,4\}}+z_{\{1,3,4\}}-z_{\{2,3,4\}}-z_{\{1,4\}} \leq 0.
$$
We now show that the above inequality is not a CG-cut for $\MPLP(G)$.
Assume, for a contradiction, that it is.
From \cref{prop: gtri} it follows that, this inequality is facet-defining for $\MP(G)$.
Therefore, it is also facet-defining for the CG-closure of $\MPLP(G)$.
This implies that the strict inequality 
$$
z_{\{1,2,4\}}+z_{\{1,3,4\}}-z_{\{2,3,4\}}-z_{\{1,4\}} < 1
$$
is valid for $\MPLP(G)$.
However, the vector $\bar z \in \MPLP(G)$ defined by $\bar z_i = \frac{1}{2}$ for $i \in \{1,2,3,4\}$, $\bar z_{\{1,2,4\}}=z_{\{1,3,4\}} = \frac{1}{2}$ and $\bar z_e = 0$ for all $e \in E \setminus \{e_1, e_3\}$ satisfies
$$
\bar z_{\{1,2,4\}}+\bar z_{\{1,3,4\}}-\bar z_{\{2,3,4\}}-\bar z_{\{1,4\}} = 1.
$$
$\diamond$
\end{example}

We would like to conclude this section by acknowledging that we are leaving several questions as directions of future research. First, an important question is how to generalize inequalities~\eqref{gtrineq} for $\alpha$-cycles of any length. The resulting inequalities will serve as the generalization of odd-cycle inequalities~\cite{Pad89}. The next step is to understand the complexity of separation over these inequalities. 
We can then define a new relaxation of the multilinear polytope by adding the proposed $\alpha$-cycle inequalities to the complete edge relaxation and characterize the class of hypergraphs for which the new relaxation coincides with the multilinear polytope. Finally, examining the computational benefits of the generalized triangle inequalities is a topic of future research.

\ifthenelse {\boolean{MPA}}
{

\bigskip
\section*{Declarations}

\noindent
\textbf{Competing interests:}
A. Del Pia is on the Editorial Board of Mathematical Programming journal. The authors have no other competing interests to declare that are relevant to the content of this article.

} {
}

\medskip
\noindent
\textbf{Funding:}
A. Del Pia is partially funded by AFOSR grant FA9550-23-1-0433 and ONR grant N00014-25-1-2490. 
A. Khajavirad is in part supported by AFOSR grant FA9550-23-1-0123 and by ONR grant N00014-25-1-2491.
Any opinions, findings, and conclusions or recommendations expressed in this material are those of the authors and do not necessarily reflect the views of the Air Force Office of Scientific Research or the Office of Naval Research.

\ifthenelse {\boolean{MPA}}
{
\bibliographystyle{spmpsci}
}
{
\bibliographystyle{plainurl}
}


\begin{thebibliography}{10}

\bibitem{BarMah86}
F.~Barahona and A.R. Mahjoub.
\newblock On the cut polytope.
\newblock {\em Mathematical Programming}, 36(2):157--173, 1986.

\bibitem{BeeFagMaiYan83}
C.~Beeri, R.~Fagin, D.~Maier, and M.~Yannakakis.
\newblock On the desirability of acyclic database schemes.
\newblock {\em Journal of the ACM}, 30:479--513, 1983.

\bibitem{BieMun18}
D.~Bienstock and G.~Mu\~noz.
\newblock {LP} formulations for polynomial optimization problems.
\newblock {\em SIAM Journal on Optimization}, 28(2):1121--1150, 2018.

\bibitem{BonGunLin18}
P.~Bonami, O.~G\"unl\"uk, and J.~Linderoth.
\newblock Globally solving nonconvex quadratic programming problems with box
  constraints via integer programming methods.
\newblock {\em Mathematical Programming Computation}, 10(3):333--382, 2018.

\bibitem{BorCraHam92}
E.~Boros, Y.~Crama, and P.L. Hammer.
\newblock Chv\'atal cuts and odd cycle inequalities in quadratic $0-1$
  optimization.
\newblock {\em SIAM Journal on Discrete Mathematics}, 5(2):163--177, 1992.

\bibitem{BorHam02}
E.~Boros and P.L. Hammer.
\newblock Pseudo-{B}oolean optimization.
\newblock {\em Discrete Applied Mathematics}, 123(1):155--225, 2002.

\bibitem{BucRin08}
C.~Buchheim and G.~Rinaldi.
\newblock Efficient reduction of polynomial zero-one optimization to the
  quadratic case.
\newblock {\em SIAM Journal on Optimization}, 18(4):1398--1413, 2008.

\bibitem{yc93}
Y.~Crama.
\newblock {Concave extensions for non-linear $\01$ maximization problems}.
\newblock {\em Mathematical Programming}, 61:53--60, 1993.

\bibitem{dPDiG21IJO}
A.~Del~Pia and S.~Di~Gregorio.
\newblock Chv\'atal rank in binary polynomial optimization.
\newblock {\em INFORMS Journal on Optimization}, 3(4):315--349, 2021.

\bibitem{dPDiG23ALG}
A.~Del~Pia and S.~Di~Gregorio.
\newblock On the complexity of binary polynomial optimization over acyclic
  hypergraphs.
\newblock {\em Algorithmica}, 85:2189--2213, 2023.

\bibitem{dPKha17MOR}
A.~Del~Pia and A.~Khajavirad.
\newblock A polyhedral study of binary polynomial programs.
\newblock {\em Mathematics of Operations Research}, 42(2):389--410, 2017.

\bibitem{dPKha18SIOPT}
A.~Del~Pia and A.~Khajavirad.
\newblock The multilinear polytope for acyclic hypergraphs.
\newblock {\em SIAM Journal on Optimization}, 28(2):1049--1076, 2018.

\bibitem{dPKha18MPA}
A.~Del~Pia and A.~Khajavirad.
\newblock On decomposability of multilinear sets.
\newblock {\em Mathematical Programming, Series A}, 170(2):387--415, 2018.

\bibitem{dPKha21MOR}
A.~Del~Pia and A.~Khajavirad.
\newblock The running intersection relaxation of the multilinear polytope.
\newblock {\em Mathematics of Operations Research}, 46(3):1008--1037, 2021.

\bibitem{dPKha23mMPA}
A.~Del~Pia and A.~Khajavirad.
\newblock A polynomial-size extended formulation for the multilinear polytope
  of beta-acyclic hypergraphs.
\newblock {\em Mathematical Programming Series A}, pages 1--33, 2023.

\bibitem{DezLau97}
M.M. Deza and M.~Laurent.
\newblock {\em Geometry of Cuts and Metrics}, volume~15 of {\em Algorithms and
  Combinatorics}.
\newblock Springer-Verlag, Berlin Heidelberg New York, 1997.

\bibitem{Dur12}
D.~Duris.
\newblock Some characterizations of $\gamma$ and $\beta$-acyclicity of
  hypergraphs.
\newblock {\em Information Processing Letters}, 112(16):617--620, 2012.

\bibitem{fagin83}
R.~Fagin.
\newblock Degrees of acyclicity for hypergraphs and relational database
  schemes.
\newblock {\em Journal of the ACM (JACM)}, 30(3):514--550, 1983.

\bibitem{GloWol74}
F.~Glover and E.~Woolsey.
\newblock Converting the 0-1 polynomial programming problem to a 0-1 linear
  program.
\newblock {\em Operations Research}, 22(1):180--182, 1974.

\bibitem{Kha22}
A.~Khajavirad.
\newblock On the strength of recursive {McCormick} relaxations for binary
  polynomial optimization.
\newblock {\em Operations Research Letters}, 51(2):146--152, 2023.

\bibitem{KhaWang25}
A.~Khajavirad and Y.~Wang.
\newblock Inference in higher order undirected graphical models and binary
  polynomial optimization.
\newblock {\em INFORMS Journal on Computing}, 2025.

\bibitem{Lau09}
M.~Laurent.
\newblock Sums of squares, moment matrices and optimization over polynomials.
\newblock In {\em Emerging Applications of Algebraic Geometry}, volume 149 of
  {\em The IMA Volumes in Mathematics and its Applications}, pages 157--270.
  Springer, 2009.

\bibitem{Pad89}
M.~Padberg.
\newblock The {B}oolean quadric polytope: Some characteristics, facets and
  relatives.
\newblock {\em Mathematical Programming}, 45(1--3):139--172, 1989.

\bibitem{Pun22book}
Abraham~P. Punnen, editor.
\newblock {\em The Quadratic Unconstrained Binary Optimization Problem: Theory,
  Algorithms and Applications}.
\newblock Springer, 2022.

\bibitem{SchBookIP}
A.~Schrijver.
\newblock {\em Theory of Linear and Integer Programming}.
\newblock Wiley, Chichester, 1986.

\bibitem{SchWal24}
E.~Schutte and M.~Walter.
\newblock Relaxation strength for multilinear optimization: Mccormick strikes
  back.
\newblock In {\em International Conference on Integer Programming and
  Combinatorial Optimization}, pages 393--404. Springer, 2024.

\bibitem{SheAda90}
H.D. Sherali and W.P. Adams.
\newblock A hierarchy of relaxations between the continuous and convex hull
  representations for zero-one programming problems.
\newblock {\em SIAM Journal of Discrete Mathematics}, 3(3):411--430, 1990.

\bibitem{wang05}
J.~Wang.
\newblock On axioms contituting the foundation of hypergraph theory.
\newblock {\em Acta Mathematicae Applicatae Sinica}, 21:495--498, 2005.

\end{thebibliography}

\end{document}